\documentclass[11pt]{article}
\usepackage{amsmath,amsthm,amssymb}
\usepackage{textcomp}
\usepackage{mathtools,bbm,xspace}
\usepackage{mathrsfs}
\usepackage{nicefrac}
\usepackage[numbers]{natbib}
\usepackage[usenames,dvipsnames]{xcolor}
\usepackage{pgfplots}
\pgfplotsset{compat=1.18}
\usepackage{ifxetex}
\usepackage{dirtytalk}
\usepackage{soul}

\usepackage[pagebackref]{hyperref}
\hypersetup{
colorlinks=true,
urlcolor=blue,
linkcolor=blue,
citecolor=OliveGreen,
}

\usepackage[capitalise,nameinlink]{cleveref}
\crefname{lemma}{Lemma}{Lemmas}
\crefname{fact}{Fact}{Facts}
\crefname{theorem}{Theorem}{Theorems}
\crefname{corollary}{Corollary}{Corollaries}
\crefname{conjecture}{Conjecture}{Conjectures}
\crefname{claim}{Claim}{Claims}
\crefname{example}{Example}{Examples}
\crefname{problem}{Problem}{Problems}
\crefname{setting}{Setting}{Settings}
\crefname{definition}{Definition}{Definitions}
\crefname{assumption}{Assumption}{Assumptions}
\crefname{subsection}{Subsection}{Subsections}
\crefname{section}{Section}{Sections}

\usepackage{tcolorbox}
\usepackage{enumerate}
\usepackage{csquotes}

\usepackage{algorithm}
\usepackage{algpseudocode}
\usepackage{tabularx}

\usepackage{todonotes}

\usepackage{bm}
\usepackage{aliascnt}
\allowdisplaybreaks

\newtheorem{theorem}{Theorem}[section]
\newtheorem*{theorem*}{Theorem}

\newtheorem*{proposition*}{Proposition}
\newaliascnt{lemma}{theorem}
\newtheorem{lemma}[lemma]{Lemma}
\aliascntresetthe{lemma}
\newtheorem*{lemma*}{Lemma}
\newtheorem{corollary}[theorem]{Corollary}
\newtheorem*{conjecture*}{Conjecture}

\newtheorem*{fact*}{Fact}

\newtheorem*{exercise*}{Exercise}

\newtheorem*{hypothesis*}{Hypothesis}

\usepackage{thm-restate}

\makeatletter
\@ifundefined{newcounteralias}{}{%
  \renewcommand\thmt@autorefsetup{%
    \@xa\def\csname\thmt@envname autorefname\@xa\endcsname
      \@xa{\thmt@thmname}%
  }%
}
\makeatother

\crefname{thm2}{Theorem}{Theorems}
\Crefname{thm2}{Theorem}{Theorems}

\theoremstyle{definition}

\newtheorem{exercise-easy}[theorem]{Exercise}
\newtheorem{exercise-med}[theorem]{Exercise}
\newtheorem{exercise-hard}[theorem]{Exercise$^\star$}
\newtheorem{claim}[theorem]{Claim}
\newtheorem*{claim*}{Claim}

\newtheorem{remark}[theorem]{Remark}
\newtheorem*{remark*}{Remark}
\newtheorem{observation}[theorem]{Observation}
\crefname{observation}{Observation}{Observations}
\Crefname{observation}{Observation}{Observations}
\newtheorem*{observation*}{Observation}

\newcommand{\Sref}[1]{\hyperref[#1]{\S\ref*{#1}}}

\newcommand\numberthis{\addtocounter{equation}{1}\tag{\theequation}}

\usepackage{parskip}
\usepackage[
letterpaper,
top=1in,
bottom=1in,
left=1in,
right=1in]{geometry}
\DeclareMathOperator*{\argmin}{arg\,min}

\DeclarePairedDelimiterX{\norm}[1]{\lVert}{\rVert}{#1}

\DeclarePairedDelimiterX{\abs}[1]{\lvert}{\rvert}{#1}
\DeclarePairedDelimiterX{\inp}[2]{\langle}{\rangle}{#1, #2}
\DeclarePairedDelimiter{\ceil}{\lceil}{\rceil}
\DeclarePairedDelimiter{\floor}{\lfloor}{\rfloor}
\DeclarePairedDelimiterX{\infdivx}[2]{(}{)}{%
  #1\;\delimsize\|\;#2%
}

\DeclareMathOperator*{\ind}{\mathbb{I}}

\DeclareMathOperator*{\E}{\mathbb{E}}
\DeclareMathOperator*{\e}{\mathbb{E}}

\DeclareMathOperator*{\Cov}{\mathrm{Cov}}
\let\Pr\relax
\DeclareMathOperator*{\Pr}{\mathbf{Pr}}
\DeclareMathOperator*{\p}{\mathbf{Pr}}

\newcommand{\lprp}[1]{\left(#1\right)}
\newcommand{\lbrb}[1]{\left\{#1\right\}}
\newcommand{\lsrs}[1]{\left[#1\right]}

\newcommand{\mc}[1]{\mathcal{#1}}
\newcommand{\mb}[1]{\mathbb{#1}}
\newcommand{\mrm}[1]{\mathrm{#1}}

\newcommand{\eps}{\varepsilon}

\newcommand{\wt}[1]{\widetilde{#1}}
\newcommand{\wh}[1]{\widehat{#1}}

\newcommand{\round}[1][t]{\left\lceil \frac{#1}{2} \right\rceil}

\newcommand{\medproj}{\mathrm{Med}_{\mathrm{proj}}}
\newcommand{\med}{\mathrm{Med}}
\newcommand{\bp}{\bar{p}}

\newcommand{\riskp}{\mathfrak R(n,p,d)}

\def\com{1}
\newcommand{\ishaq}[1]{
    \if\com1
        \todo[inline,color=blue!30]{\small\textbf{Ishaq:} #1}
    \else
    \fi
}
\newcommand{\kasper}[1]
{
    \if\com1
        \todo[inline,color=green!30]{\small\textbf{Kasper:} #1}
    \else
    \fi
}
\newcommand{\mikael}[1]
{
    \if\com1
        \todo[inline,color=red!30]{\small\textbf{Mikael:} #1}
    \else
    \fi
}

\newcommand{\nikita}[1]
{
    \if\com1
        \todo[inline,color=yellow!30]{\small\textbf{Nikita:} #1}
    \else
    \fi
}

\newcommand{\yeshwanth}[1]
{
    \if\com1
        \todo[inline,color=orange!30]{\small\textbf{Yeshwanth:} #1}
    \else
    \fi
}

\crefformat{equation}{(#2#1#3)}

\makeatletter
\newcommand{\alphacmd@factory}[1]{}
\newcounter{alphacmdcounter}
\newcommand{\GenerateAlphabetCmds}[2]{%
    \renewcommand{\alphacmd@factory}[1]{%
        \expandafter\providecommand\csname #1##1\endcsname{{#2{##1}}}%
    }
    \setcounter{alphacmdcounter}{0}
    \loop
        \stepcounter{alphacmdcounter}
        \edef\alphacmd@ID{\@Alph\c@alphacmdcounter}
        \expandafter\alphacmd@factory\alphacmd@ID
    \ifnum\thealphacmdcounter<26
    \repeat
}
\newcommand{\GenerateAlphabetCmdsLower}[2]{%
    \renewcommand{\alphacmd@factory}[1]{%
        \expandafter\providecommand\csname #1##1\endcsname{{#2{##1}}}%
    }
    \setcounter{alphacmdcounter}{0}
    \loop
        \stepcounter{alphacmdcounter}
        \edef\alphacmd@ID{\@alph\c@alphacmdcounter}
        \expandafter\alphacmd@factory\alphacmd@ID
    \ifnum\thealphacmdcounter<26
    \repeat
}
\makeatother

\GenerateAlphabetCmds{c}{\mathcal}
\GenerateAlphabetCmdsLower{c}{\mathcal}

\GenerateAlphabetCmds{r}{\mathbf}
\GenerateAlphabetCmdsLower{r}{\mathbf}

\usepackage{graphicx}
\title{
Phase Transitions in Heavy-Tailed Mean Estimation under $\ell_p$ Norms
}

 \author{Ishaq Aden-Ali\thanks{Work done at EECS Department, UC Berkeley. Email: {adenali@berkeley.edu}} \and Yeshwanth Cherapanamjeri \thanks{College of Computing and Data Science, Nanyang Technological University, Email: yesh.c@ntu.edu.sg}\and Mikael Møller Høgsgaard\thanks{Department of Statistics, University of Oxford and Computer Science Department, Aarhus University. Email: hogsgaard@cs.au.dk} \and Kasper Green Larsen\thanks{Computer Science Department, Aarhus University. Email:  larsen@cs.au.dk} \and Nikita Zhivotovskiy\thanks{Department of Statistics, UC Berkeley. Email: zhivotovskiy@berkeley.edu}}
\date{May 6, 2026}
\begin{document}

\maketitle

\begin{abstract}
We study the problem of estimating the mean of a random vector in
$\mathbb R^d$ under $\ell_p$ norms, assuming only the existence of a covariance
matrix. In the Euclidean norm, the sample mean is already optimal in expectation,
and robust estimators are usually needed only to obtain high-probability tails.
Our results show that this picture changes for $\ell_p$ norms with $p>2$:
robustification is already necessary even for the expected risk. We prove that
the minimax risk exhibits a sharp phase transition, controlled by the interaction
between the sample size $n$, the dimension $d$, and the norm parameter $p$. In
the large-sample regime, the classical Gaussian rate remains achievable. In the
complementary regime, heavy-tailed distributions create a genuinely non-Gaussian
obstruction. Furthermore we show that the simple estimator: coordinate-wise median of means whose number of buckets chosen according to the geometry of the loss, optimally tracks this phase transition in its expected risk minimax. We also show that natural
projection-based median-of-means estimators for general norms, despite achieving
the optimal confidence trade-off relative to the sample-mean benchmark of Lugosi
and Mendelson (\emph{Probab. Theory Relat. Fields}, 2019), can be polynomially
suboptimal for $\ell_p$ norms even in the constant-probability regime.
\end{abstract}

\section{Introduction}
Mean estimation, or location estimation, is one of the foundational subroutines in machine learning and statistics. A modern line of work focuses on statistical estimators that are robust to heavy-tailed distributions.
Starting with the systematic study pioneered by Catoni \cite{catoni2012challenging}, this line of work seeks to obtain non-asymptotic \emph{sub-Gaussian} rates for estimating the means of random variables, usually using estimators such as the median of means, which appeared independently in
\cite{NemirovskyYudin1983ProblemComplexity,JerrumValiantVazirani1986RandomGeneration,AlonMatiasSzegedy1999FrequencyMoments}, Catoni's estimator \cite{catoni2012challenging}, or trimmed-mean estimators \cite{LugosiMendelson2021TrimmedMean}. As we recall later, more involved versions of these estimators also work well in high-dimensional regimes, where the primary focus has traditionally been on the Euclidean norm.

However, the Euclidean norm is arguably the most benign norm to work with. Indeed, given $n$ independent vectors $X_1,\ldots,X_n$ in $\mathbb{R}^d$ with common mean $\mu$ and identity covariance $I_d$, the sample mean is essentially optimal, up to universal constant factors for all $n$ and $d$:
\begin{equation}
\label{eq:simplesamplemean}
\E\left\|\frac{1}{n}\sum_{i=1}^n X_i-\mu\right\|_2
\le \sqrt{\frac{d}{n}},
\end{equation}
which follows from Cauchy--Schwarz and does not require the existence of any moments above the second moment, while the lower bound follows from Anderson's lemma \cite{Anderson1955SymmetricUnimodal} as discussed below. Assuming only the existence of the covariance, and no higher moments, the main focus has therefore been on obtaining sub-Gaussian estimators, namely estimators with optimal sub-Gaussian type tails, where the sample mean is clearly the wrong object
\cite{Minsker2015GeometricMedian,LugosiMendelson2019SubGaussian,CatoniGiulini2017DimensionFreePACBayesian,LugosiMendelson2021TrimmedMean}.

When going beyond the Euclidean norm, while still keeping only the bounded covariance assumption, even the in expectation story of \eqref{eq:simplesamplemean} becomes more delicate especially in the high-dimensional regime where the dimension $d$ is larger than the sample size $n$. In non-Hilbertian geometries, the vanilla sample mean is no longer necessarily a good estimator.\footnote{As an example, for matrix valued observations the natural benchmark is often the operator norm. In this setting, heavy-tailed matrices are known to require additional truncation steps even for in expectation risk bounds, before one can ask for sub-Gaussian tails; see, for example, \cite{AbdallaZhivotovskiy2026Covariance}.}
Our key contribution is to show that the Euclidean picture is special even for the simple mean estimation problem in $\mathbb{R}^d$. For $\ell_p$ norms with $p>2$, robustification can be necessary even for \emph{in expectation} bounds of the form \eqref{eq:simplesamplemean}. This phenomenon does not appear in the regime $1\le p\le 2$, where the sample mean remains optimal up to universal constants, but it becomes crucial once the geometry of the loss is sufficiently far from $\ell_2$.

\begin{theorem}[Phase transition for location estimation in $\ell_p$, informal]
\label{thm:mainsimple}
Under only the covariance assumption $\operatorname{Cov}_P(X)=I_d$, the optimal expected error for estimating the mean in $\ell_p$ satisfies
\[
\inf_{\widehat{\mu}}
\sup_{\substack{ \mu \in \mathbb{R}^d \\ P:\, \E_P [X]=\mu\\ \operatorname{Cov}_P(X)=I_d}}
\E\|\widehat{\mu}-\mu\|_p
\asymp
\begin{cases}
 d^{\frac{1}{2\lceil n/2\rceil}},
 &  \bar p > 2\lceil n/2\rceil,\\[1ex]
 d^{\frac{1}{\bar p}}\sqrt{\frac{\bar p}{n}},
 &  \bar p \le 2\lceil n/2\rceil,
\end{cases}
\]
where $ \hat{\mu} $ ranges over all estimators based on the $ n $ samples and
\[
\bar p = \min \{ p,
\, \max\{\log d,\,  1\} \}.
\]
\end{theorem}
Thus, the minimax risk undergoes a sharp phase transition. In the large sample regime, the rate has the familiar Gaussian form, corresponding to the Gaussian width of the dual ball inducing the norm, as we discuss below:
\[
d^{1/\bar p}\sqrt{\frac{\bar p}{n}}.
\]
In contrast, in the small sample or large $p$ regime, the effective exponent is capped by the number of samples, and the optimal rate becomes
\[
d^{\frac{1}{2\lceil n/2\rceil}}.
\]
We show that a estimator achieving this rate is a simple coordinate-wise median of means estimator, with the choice of number of buckets depending on the geometry of the loss. Recall that median of means first splits the sample into nearly equal buckets, computes the average inside each bucket, and then takes a median of these bucket averages. In the coordinate-wise version used here, if $\bar X^{(1)},\ldots,\bar X^{(t)}$ denote the nonempty bucket averages, the estimator is defined by
\[
\widehat \mu_j
=
\operatorname{Median}\left(\bigl\{
\bar X^{(1)}_j,\ldots,\bar X^{(t)}_j
\bigr\}\right),
\qquad j=1,\ldots,d.
\]
The number of buckets is dictated by the geometry of the loss, namely by the effective exponent $\bar p$. More precisely, the estimator uses
\[
t(\bar p)=2\left\lceil \frac{\bar p}{2}\right\rceil-1
\]
buckets, with the median taken over the nonempty buckets. When $1\le p<2$, we have $t(\bar p)=1$, so the estimator is just the ordinary sample mean. The same happens at $p=2$, recovering exactly the Euclidean behavior in \eqref{eq:simplesamplemean}. When $p>2$ and $\bar p>2$, the estimator becomes a genuine coordinate-wise median of means estimator, with the number of buckets increasing with $\bar p$. Finally, in the extreme regime where $t(\bar p)>n$, every nonempty bucket contains a single observation, and the procedure reduces to the coordinate-wise median of the sample points. The technical difficulty is to show that these regimes match sharply and together give the precise minimax rate above.

\paragraph{Phase transition and the existing general norm results.}
Our \emph{in expectation} analysis complements the existing high probability theory of robust mean estimation with respect to general norms. Sub-Gaussian mean estimators for general norms were studied, in different forms, in
\cite{LugosiMendelson2019GeneralNorms,Minsker2025UniformBounds,DepersinLecue2021PseudoNorms}. The original benchmark in \cite{LugosiMendelson2019GeneralNorms} is distribution dependent and is centered around robustly reproducing the performance of the empirical mean, together with the correct one dimensional sub-Gaussian deviation term. Informally, in our $\ell_p$ setup this benchmark contains a term of the form
\[
\E\left\|
\frac{1}{n}\sum_{i=1}^n (X_i-\mu)
\right\|_p ,
\]
or its symmetrized analogue. This is the natural target if the goal is to obtain optimal high probability tails around the performance of the sample mean. However, for the problem considered in this paper, this benchmark can be too large: as follows from our results, for $p>2$, the sample mean itself is already suboptimal even in expectation.\footnote{For example, let $X=\eps\sqrt{d}\,e_J$, where $\eps$ is Rademacher and $J$ is uniform on $[d]$. Then $\E X=0$ and $\Cov(X)=I_d$. If $n\leq c\sqrt{d}$, the sampled basis vectors are distinct with constant probability, so
$
\E\left\|n^{-1}\sum_{i=1}^n X_i\right\|_p
\gtrsim \frac{\sqrt{d}}{n^{1-1/p}}.
$
When $2<p\leq\log d$ and $p\leq2\lceil n/2\rceil$, \cref{thm:mainsimple} gives the rate $d^{1/p}\sqrt{p/n}$. Their ratio is
$\asymp (d/n)^{1/2-1/p}/\sqrt{p}$, which diverges for fixed $p>2$ as $d/n\to\infty$. Thus, the sample mean is suboptimal when $n\leq c\sqrt{d}$ and $2<p\leq\min\{\log d,2\lceil n/2\rceil\}$.}

A related but different point is emphasized in \cite{DepersinLecue2021PseudoNorms}. If one takes Gaussian mean estimation as the benchmark, then the relevant complexity is the Gaussian width of the dual ball inducing the norm. In the case of the $\ell_p$ norm and identity covariance, this gives
\[
\frac{\E\|G\|_p}{\sqrt n}
\asymp
d^{1/\bar p}\sqrt{\frac{\bar p}{n}},
\]
where $G\sim N(0,I_d)$ and $\bar p=\min \{ p,
\, \max\{\log d,\,  1\} \}$. Thus, the second branch of our theorem,
\[
d^{1/\bar p}\sqrt{\frac{\bar p}{n}},
\]
coincides exactly with the Gaussian width benchmark. In this regime, our estimator achieves the same rate that the sample mean would have for Gaussian data; this may be called the sub-Gaussian regime. The genuinely new phenomenon provided by our analysis is the first branch,
\[
d^{\frac{1}{2\round[n]}}.
\]
This term does not appear in the Gaussian mean estimation setup. It reflects a finite sample heavy-tailed regime under only a covariance assumption. Morally, the Gaussian width rate is correct once the number of samples is large enough relative to the effective exponent $\bar p$. When this condition fails, the minimax risk enters a different phase.

\paragraph{Phase transition and the Gaussian location estimation problem.}
The Gaussian width branch of Theorem \ref{thm:mainsimple} is also closely connected to the classical minimax theory of Gaussian location estimation. In the Gaussian location model, Anderson's lemma implies an exact minimax identity for every bowl-shaped loss; see the original inequality of Anderson
\cite{Anderson1955SymmetricUnimodal} and the modern presentation in
\cite[Section~28.6]{PolyanskiyWu2025InformationTheory}. In particular, for
$X_1,\ldots,X_n\sim N(\mu,I_d)$, Anderson's lemma implies
\[
\inf_{\widehat\mu}\sup_{\mu\in\mathbb R^d}
\E_\mu \|\widehat\mu-\mu\|_p
=
\frac{\E\|G\|_p}{\sqrt n}
\asymp
d^{1/\bar p}\sqrt{\frac{\bar p}{n}},
\qquad
G\sim N(0,I_d),
\quad
\bar p=\min\{ p ,\, \max\{\log d,\, 1\} \}.
\]
Thus, as discussed above, the second branch of our theorem is exactly the classical Gaussian location obstruction. Our result shows that this Gaussian answer is not the whole story under only a covariance assumption. For the much larger class of all distributions with covariance $I_d$, the Gaussian benchmark remains sharp only when
\[
\bar p\le 2\round[n] .
\]
When this condition fails, heavy-tailed distributions create an additional obstruction, and the minimax risk becomes
$
d^{\frac{1}{2\round[n]}}.
$
In this sense, Theorem \ref{thm:mainsimple} may be viewed as a heavy-tailed, covariance-only analogue of the classical Gaussian location minimax identity.

\paragraph{High probability implications and limitations of projection-based estimators.}
As of now, the full high probability minimax picture for $\ell_p$ norms remains unclear. In
particular, the aforementioned general norm results do not directly provide the matching
high probability behavior for $\ell_p$ norms under only the covariance assumption. Our theorem in expectation already shows that the sample mean can be suboptimal. On the negative side, we further
show that natural projection-based median of means estimators, of the type used in the
general norm theory
\cite{LugosiMendelson2019GeneralNorms,Minsker2025UniformBounds,DepersinLecue2021PseudoNorms,OliveiraResende2025TrimmedSampleMeans},
can also fail to achieve the correct $\ell_p$ behavior.

Informally, a projection-based median of means estimator works as follows. It searches for a point
whose projection in every dual direction agrees with the median of the projected bucket means:
\[
\widehat \mu_{\rm proj}
\in
\arg\min_x
\sup_{\|v\|_q\le 1}
\left|
\langle v,x\rangle
-
\operatorname{Median}
 \left(\bigl\{
\langle v,\bar X^{(1)}\rangle,\ldots,
\langle v,\bar X^{(t)}\rangle
\bigr\}\right)
\right|,
\qquad
\frac1p+\frac1q=1.
\]
This construction is natural as a generic norm-agnostic procedure, but our result shows that it can be too crude even for the classical family of $\ell_p$ norms once $p>2$.
\begin{restatable}[Suboptimality of Projection-based median-of-means]{theorem}{ProjSuboptimal}
    \label{thm:proj_sub_optimal}
    Let $p > 2$. There exist absolute constants $C, c > 0$ such that the following hold. Let $d, n \in \mathbb{N}$ be such that $n \leq c \sqrt{d}$ and $d > C$. Then, there exists a distribution $D$ with mean $0$ and covariance $I_d$ such that the projection-based median-of-means estimator with any bucket size $m$ satisfying $n / m \geq 3$, $\wh{\mu}_{\mrm{Proj}}$, satisfies:
    \begin{equation*}
        \Pr_{\bm{X} \thicksim D^{\otimes n}} \left[\norm{\wh{\mu}_\mrm{Proj} \lprp{\bm{X}}}_p \geq \frac{\sqrt{d}}{2 n^{1 - 1 / p}}\right] \geq 0.9.
    \end{equation*}
\end{restatable}

This should be compared with the coordinate-wise estimator analyzed above. In the Gaussian width
regime, the latter achieves a much smaller rate. For example, for $p=\infty$ and
$\log d\lesssim n\ll \sqrt d$, the coordinate-wise estimator has error of order, with constant
probability,
\[
\sqrt{\frac{\log d}{n}},
\]
while the projection-based median of means lower bound is of order, again with constant
probability,
\[
\frac{\sqrt d}{n}.
\]
Thus, in this regime the projection-based estimator pays a polynomial price in the dimension,
while the coordinate-wise estimator has only logarithmic dependence. This shows that, for
structured norms such as $\ell_p$, one needs to go beyond projection-based estimators for general
norms.

\paragraph{A high probability result for the $\ell_\infty$ norm.}
Although
our main focus is on the results in expectation, as a complementary result we give a sharp high probability analysis for the $\ell_\infty$ norm. This result is the
high probability analogue of the phase transition phenomenon above, and refines even the
constant probability regimes implied by our analysis in expectation.

\begin{theorem}[$\ell_\infty$ high probability phase transition, informal]
\label{thm:inf_high_prob_phase_transition}
Under only the covariance assumption $\operatorname{Cov}_P(X)= I_d$, the optimal
$\ell_\infty$ error at confidence level $1-\delta$ is, up to universal constants,
\[
\begin{cases}
 \left(\dfrac d\delta\right)^{\frac{1}{2\round[n]}},
&
n \lesssim \log(d/\delta),
\\[2ex]

\sqrt{\dfrac{\log(d/\delta)}{n}},
&
n \gtrsim \log(d/\delta).
\end{cases}
\]
The upper bound is achieved by the coordinate-wise median of means estimator with
$t=2\lceil \log(d/\delta)\rceil$ buckets, and the lower bound holds for every estimator.
\end{theorem}

As before, the first regime is genuinely non-Gaussian: the error is polynomial in $d/\delta$,
although the exponent decreases with the sample size. Once $n$ is of order $\log(d/\delta)$, the
usual sub-Gaussian rate takes over. Thus, even for the classical $\ell_\infty$ norm, the
Gaussian width or sub-Gaussian benchmark captures only the large sample side of the theory,
while under only a covariance assumption there is an additional finite sample heavy-tailed phase.
This also highlights a limitation of projection-based estimators: in the regimes where the
coordinate-wise estimator attains the sharp $\ell_\infty$ rate, projection-based median of means
estimators can be polynomially suboptimal, as shown by \cref{thm:proj_sub_optimal}.

\paragraph{On computation and coordinate-wise medians.}
A significant line of work on heavy tailed mean estimation in the Euclidean norm studies polynomial time estimators with sub-Gaussian guarantees under finite covariance assumptions; see, for example,
\cite{Hopkins2020SubGaussianPolynomialTime,CherapanamjeriFlammarionBartlett2019FastMean,HopkinsLiZhang2020RegretMean,DiakonikolasKanePensia2020Stability}. In the main regimes in expectation studied here, computation is essentially trivial: the optimal estimator is coordinate-wise and reduces to bucket averages and one dimensional medians. At the same time, coordinate-wise median of means is not universally optimal. In the high probability Euclidean problem, no single distribution independent choice of the number of buckets is optimal for all distributions with identity covariance.
We prove the following result.

\begin{restatable}[Suboptimality of coordinate-wise median-of-means]{theorem}{MoMSuboptimal}
    \label{thm:mom_l2_lb_main}
There are universal constants $c, C \ge 0$ such that the following holds.  For any $\delta \in (0, c)$ and $ n, d, t \in\mathbb{N}$ with $ t\leq n $, there exists a distribution $P$ over $\mathbb{R}^{d} $ with mean $0$ and covariance $I_d$ such that if~\cref{alg:MOM} is ran with $t$ buckets and sample $X_1, \dots, X_n \sim P$ as input, it returns an estimate $\hat\mu$ satisfying
\[
\Pr_{X_1, \dots, X_n \sim P}\left[ \|\hat \mu \|_2 \ge C \sqrt{\frac{d\log\left(1/\delta\right)}{n\log\log\left(1/\delta\right)}} \right] > \delta.
\]
\end{restatable}

This is expected, and probably folklore, but we include it for completeness in \cref{sec:suboptimalcoordwise}. The complete high probability theory for general norms therefore remains substantially more delicate than the picture in expectation developed here.\footnote{Recent generic chaining based approaches are particularly relevant here; see
\cite[Theorem~5.1]{AbdallaMendelson2026CovarianceDirectionDependent}
and \cite{BartlMendelson2026UniformMeanGenericChaining}. We note that these works do not provide a sharp high probability minimax theory for $\ell_p$ norms under the covariance only assumption.}

\section{Notation and Organization}
We begin with some notation. We use $\log$ to denote the natural logarithm. For a real number $a$, we will use $\floor{a}$ and $\ceil{a}$ to denote its floor and ceiling respectively. For a positive integer $a$, we will often use the fact that $\lfloor \frac{a + 1}{2} \rfloor$ is a short form for $\frac{a+1}{2}$  when $a$ is odd and $\frac{a}{2}$ when $a$ is even. For a multiset of real numbers, $\{a_{i}\}_{i = 1}^n$, we use $a_{(1)} \geq a_{(2)} \dots \geq a_{(n)}$ to denote the decreasing rearrangement of its elements with the median being $a_{\floor{(n + 1) / 2}}$ if $n$ is odd or $(a_{(n / 2)} + a_{(n / 2 + 1)}) / 2$ when $n$ is even. For positive integer $a$ we define $[a] = \{1, 2, \dots, a\}$.
We use $c, C$ for universal constants, and their value may change from line to line.
For $x \in \mb{R}^d$ and $p \in [1, \infty)$, we use $\norm{x}_p$ to denote the $\ell_p$ norm of $x$, with $q$ typically used to denote the dual norm exponent satisfying $\frac{1}{p} + \frac{1}{q} = 1$. We use $e_1, \dots, e_d \in \mb{R}^d$ to denote the standard basis vectors over $\mb{R}^d$. We use $I_d$ for identity matrix whose dimension will be clear from context.

For a distribution $P$ over $\mc{X}$, event $A: \mc{X} \to \lbrb{0, 1}$ and function $f: \mc{X} \to \mb{R}$, we will use $\Pr_{X \thicksim P} [ A(X) ]$, and $\E_{X \thicksim P} [f(X)]$ to be the probability of event $A$ and the expectation of $f$ when $X$ is drawn from $P$ respectively. When a random variable or its distribution is clear from context, we often drop the subscript and simply write $\Pr_P$ ($\E_P$) or $\Pr$ ($\E$). Furthermore, for a natural number $n$, we use $P^{\otimes n}$ to denote the product distribution over $n$ i.i.d.\ samples from $P$. For a distribution $P$ over $\mb{R}^d$ with mean $\mu$, we define the covariance matrix as $\Cov_{P}(X) = \E_{X \sim P} [XX^\top]- \mu\mu^\top$.

For quantities $a, b$, we write $a \lesssim b$ ($a \gtrsim b$) if there exists a universal constant $C > 0$ ($c > 0$) such that $a \leq C b$ ($a \geq cb$) uniformly over the relevant parameter ranges which will be clear from context. We write $a = O(b)$ ($a = \Omega(b)$) if $a \lesssim b$ ($a \gtrsim b$) asymptotically in the relevant parameters and $a = \Theta (b)$ if both hold. Finally, $a = o(b)$ ($a = \omega(b)$) if $\frac{a}{b}$ ($\frac{b}{a}$) goes to $0$ asymptotically.

\paragraph{Organization:} The remainder of the paper is organized as follows: \cref{thm:mainsimple} is formalized and proved in \cref{sec:opt_in_exp}, \cref{thm:inf_high_prob_phase_transition} in \cref{sec:optimal_high_prob_ell_inf}, \cref{thm:proj_sub_optimal} in \cref{ssec:app_sub_opt_proj_median},
and~\cref{thm:mom_l2_lb_main} in \cref{sec:suboptimalcoordwise}.

\section{Optimal In-expectation Guarantees}
\label{sec:opt_in_exp}

In this section we prove sharp bounds for the mini-max risk of heavy-tailed isotropic mean estimation.
The minimax risk of sample size $n$, dimension $d$ and norm $\| \cdot \| : \mb{R}^d \to \mb{R}_{\ge 0}$ is defined as follows:
\begin{equation}\label{eq:minimax-risk}
    \mathfrak R(n,d, \| \cdot \|) =
\inf_{\cA:(\mathbb R^d)^n\to\mathbb R^d}
\sup_{\substack{\mu \in \mb{R}^d \\ P:\E_P[X]=\mu,\ \Cov_P(X)= I_d}}
\E_{X_1,\dots,X_n\sim P}\bigl[\|\cA(X_1,\dots,X_n)-\mu\|\bigr].
\end{equation}
In words, this is the smallest expected error when estimating the mean, measured in the norm $\|\cdot\|$, that an estimator can guarantee based on $n$ samples, uniformly over all $d$-dimensional distributions with identity covariance. For $\ell_p$ norms, we will simply write $  \mathfrak R(n,d,p) = \mathfrak R(n,d, \| \cdot \|_p) $.

To state our main result, it will be useful to setup some notation that will be used throughout this section.
We will make frequent use of the truncated version of $p \in [1, \infty]$ defined as
\begin{equation}\label{eq:bp}
\bp = \min\{ p ,\, \max\{\log d,\, 1\} \}.
\end{equation}
The truncation at $ \log d $ is because for any $ p,q $ such that $ \log d\leq p \leq q $, we have that $ e^{-1} \|\cdot\|_{p} \leq \|\cdot\|_{q} \leq \|\cdot\|_{p} $, i.e.\ the norms are equivalent up to a constant factor. Thus, the minimax risk with respect to $ \ell_{p} $ and $ \ell_{q} $ are equivalent up to a constant factor so it suffices to consider the range $ p \in [1, \log d]$.\footnote{The truncation with $ 1 $ is just to avoid technicalities when $ d\leq 2$.}

We can now state our main result.
\begin{theorem}\label{thm:expected_minimax_risk}
For any $n, d\in\mathbb{N}$ and $p \in [1, \infty]$ we have that
\[
\riskp \asymp
\begin{cases}
 d^{\frac{1}{2\lceil n/2\rceil}}, &  \bp > 2\round[n],\\[1ex]
d^{\frac{1}{\bp}}\sqrt{\frac{\bp}{n}}, &  \bp \le 2\round[n].
\end{cases}
\]
Equivalently, if we define $s=\min\{ \bp ,\, 2\round[n]\}$, then
\[
\riskp \asymp d^{\frac{1}{s}}\sqrt{\frac{s}{n}}.
\]
Moreover, this rate is achieved by the coordinate-wise median-of-means estimator.
\end{theorem}

\begin{remark}
Observe that a phase transition occurs exactly at $\bp=2\round[n]$, where the small sample size regime $2\round[n] < \bar{p}$ rate $d^{\frac{1}{2\round[n]}}$ changes to the large sample size regime $  \bar{p} \leq 2\round[n] $ rate $d^{\frac{1}{\bp}}\sqrt{\frac{\bp}{n}}$.
\end{remark}

\subsection{Upper bound}

\begin{algorithm}[t]
\caption{Coordinate-wise Median-of-Means}\label{alg:MOM}
\begin{algorithmic}[1]
\Require Number of buckets $t\leq n$, samples $X_1,\dots,X_n\in\mathbb{R}^d$
\State Partition $[n]$ into $t$ size-balanced buckets $ B_{1},\ldots,B_{t} $ by setting $B_i = \{m\in[n] : m\equiv i \pmod t\}$
\For{$i=1$ to $t$}
    \State $\hat\mu^{(i)} \gets \frac{1}{|B_i|}\sum_{m\in B_i}X_m$
\EndFor
\State \Return $\hat\mu \gets
\operatorname{Median}
\bigl(\hat\mu^{(1)},\dots,\hat\mu^{(t)}\bigr)$ \Comment{coordinate-wise median}
\end{algorithmic}
\end{algorithm}

In this section we prove upper bounds on the expected error of the coordinate-wise median-of-means estimator (\cref{alg:MOM}) when used to estimate the mean with respect to $\ell_p$ norms for $p \in [1 , \infty]$.
We will make frequent use of the function $t(\bp)$:
\begin{equation}\label{eq:t(p)}
  t(\bp)=2\left\lceil \frac{\bp}{2}\right\rceil-1.
\end{equation}
We will use $t(\bp)$ as the number of buckets in coordinate-wise median-of-means estimator~\cref{alg:MOM} when estimating the mean with respect to the $\ell_p$ norm.
Notice that in the case that $p \in [1, 2]$, the number of buckets is just $1$ and \cref{alg:MOM} reduces to the usual sample mean that is optimal in-expectation for $p \in [1, 2]$.

We now state our upper bound on the performance of \cref{alg:MOM} for estimating the mean of a random vector in $ \mathbb{R}^{d} $ under $ \ell_{p} $ norm for any $ p\geq 1. $

\begin{theorem}\label{thm:medianofmeans}
Let $X_{1},\ldots,X_{n} \in \mb{R}^d$ be i.i.d.\ random vectors with finite mean $\mu$ and covariance $\Sigma$.
Fix $ p\geq 1$.
Given the number of buckets $t(\bar{p})$ and the sample $X_1, \ldots, X_n $ as input, \cref{alg:MOM} returns an estimate $\hat\mu$ satisfying
\[
\E\left[\|\hat\mu-\mu\|_p\right]
\le
\begin{cases}
\displaystyle
2\left(\sum_{j=1}^d \Sigma_{jj}^{\round[n]}\right)^{\frac{1}{2\round[n]}},
&\bp > 2\round[n],\\[2ex]
\displaystyle
2\,\frac{\left(\sum_{j=1}^d \Sigma_{jj}^{\bp/2}\right)^{\frac{1}{\bp}}}{\sqrt{\left\lfloor n/t(\bp)\right\rfloor}},
& \bp \le 2\round[n].
\end{cases}
\]
Moreover, if we define $s=\min\{ \bp ,\, 2\round[n]\}$ then
\[
\E\left[\|\hat\mu-\mu\|_p\right]
\le
4
\left(\sum_{j=1}^d \Sigma_{jj}^{s/2}\right)^{\frac{1}{s}}
\sqrt{\frac{s}{n}}.
\]
\end{theorem}

To prove \cref{thm:medianofmeans}, we will need the following moment bound on the deviation of the empirical median from the population mean which is a slight generalization of a result due to \cite{LarsenPT21}.
For completeness, we include a proof in \cref{app:omitted}.

\begin{restatable}{thm2}{momentsmedian}\label{momentsmedian}
Let $ X_1,\ldots,X_t $ be $t$ independent random variables with a common mean $\mu$ and let $ M = M(X_1, \dots, X_t)$ denote their median.
For all $ q\geq 0 $ it holds that
\begin{align*}
 \e\left[|M - \mu |^{q \round }\right]\leq \binom{t}{\round }\max_{i \in [t]}\e\left[|X_i-\mu|^{q}\right]^{\round}.
\end{align*}
\end{restatable}

\begin{proof}[Proof of \cref{thm:medianofmeans}]
We consider two cases.
First, let us consider the case that $\bp \le 2\round[n]$.
Notice that this condition is equivalent to the condition that $n \ge t(\bp)$, which will be more convenient to use in what follows.
Indeed, we have that
\begin{equation}\label{eq:equivalence}
t(\bp) \le n \iff \round[\bp] \le \left\lfloor \frac{n+1}{2}  \right\rfloor \iff  \round[\bp] \le \round[n] \iff \frac{\bp}{2} \le \round[n] \iff \bp \le 2\round[n].
\end{equation}
Now, using the fact that $\| \cdot \|_p \le \| \cdot \|_{\bar p}$ together with Jensen's inequality yields
\begin{equation}\label{eq:first_jensen}
\E \left[\| \hat\mu - \mu \|_p \right] \le \E \left[\| \hat\mu - \mu \|_{\bp} \right] \le \E \left[\| \hat\mu - \mu \|_{\bar p}^{\bp }\right]^{\frac{1}{\bp}} = \left( \sum_{j=1}^d \E\left[|\bar{\mu}_j - \mu_j|^{\bp} \right]\right)^{\frac{1}{\bp}}.
\end{equation}
It now suffices to bound  the terms $ \E \left[|\hat\mu_j - \mu_j|^{\bp}\right]$.
Notice that $\round[t(\bp)] = \lfloor \frac{t(\bp)+1}{2} \rfloor = \lceil \frac{\bp}{2}\rceil \ge 1$, so there is some $q \in [1 ,2]$ such that $q\round[t(\bp)]=\bp$. Thus, we can invoke the moment bound of~\cref{momentsmedian} with $t = t(\bp)$ and $q \in [1, 2]$ as above together with Jensen's inequality to get
\begin{align*}
\E \left[|\hat\mu_j - \mu_j|^{\bp}\right] &\le \binom{t(\bp)}{\round[t(\bp)]} \max_{i \in [t(\bp)]} \E \left[|\hat\mu^{(i)}_j - \mu_j|^q \right]^{\round[t(\bp)]} \\
&\le   \binom{t(\bp)}{\round[t(\bp)]}  \max_{i \in [t(\bp)]} \E \left[|\hat\mu^{(i)}_j - \mu_j|^2 \right]^{\frac{q}{2}\round[t(\bp)]}\\
&\le \binom{t(\bp)}{\round[t(\bp)]} \left(\frac{\Sigma_{jj}}{\floor{\frac{n}{t(\bp)}}}\right)^{\frac{q}{2}\round[t(\bp)]} = \binom{t(\bp)}{\floor{\frac{t(\bp)+1}{2}}} \left(\frac{\Sigma_{jj}}{\floor{\frac{n}{t(\bp)}}}\right)^{\frac{\bp}{2}},
\end{align*}
where the last inequality follows from the fact that each bucket gets at least $\floor{n/t(\bp)}$ samples and that the variance of the $j$th coordinate is $\Sigma_{jj}$.
Plugging the above back into \cref{eq:first_jensen} yields
\[
\E \left[\| \hat\mu - \mu \|_p \right] \le \binom{t(\bar{p})}{\round[t(\bp)]}^{\frac{1}{\bar{p}}}\frac{\left( \sum_{j=1}^d\left(\Sigma_{jj}\right)^{\frac{\bp}{2}}\right)^{\frac{1}{\bp}}}{\sqrt{\floor{n/t(\bp)}}} \le 2\frac{\left( \sum_{j=1}^d\left(\Sigma_{jj}\right)^{\frac{\bp}{2}}\right)^{\frac{1}{\bp}}}{\sqrt{\floor{n/t(\bp)}}},
\]
where the final inequality follows from the simple bound
\[
\binom{t(\bp)}{\round[t(\bp)]}= \binom{2\ceil{\bp/2}-1}{\ceil{\bp/2}} = \binom{2\ceil{\bp/2}-1}{\ceil{\bp/2}-1} \le \sum_{k=0}^{\ceil{\bp/2}-1} \binom{2\ceil{\bp/2}-1}{k} = 2^{2\ceil{\frac{\bp}{2}}-2} \le 2^{\bp}.
\]
Now consider the case that $\bp > 2\round[n]$, which is equivalent to $n < t(\bp)$ by \cref{eq:equivalence}. In this case, every non-empty bucket has a single sample and the estimator is the coordinate-wise median.
Since $n<t(\bp)=2\lceil \bp/2\rceil-1$, a simple case analysis on the parity of $n$ shows that $p > 2 \ceil{\frac{n}{2}} \eqqcolon r(n)$.
Thus, using the fact that $\| \cdot \|_p \le \| \cdot \|_{2\ceil{\frac{n}{2}}}$ and Jensen's inequality yields
\begin{equation}\label{eq:first_jensen_small_n}
\E \left[\| \hat\mu - \mu \|_p \right] \le \E \left[\| \hat\mu - \mu \|_{2\ceil{\frac{n}{2}}} \right] \le \E \left[\| \hat\mu - \mu \|_{r(n)}^{r(n) }\right]^{\frac{1}{r(n)}} = \left( \sum_{j=1}^d \E\left[|\bar{\mu}_j - \mu_j|^{r(n)} \right]\right)^{\frac{1}{r(n)}}.
\end{equation}
So, it suffices to bound the terms $ \E \left[|\hat\mu_j - \mu_j|^{r(n)}\right]$.
As in the previous case, we can invoke \cref{momentsmedian} with $q=2$ and $t = n$ to get the bound
\begin{align*}
\E \left[|\hat\mu_j - \mu_j|^{r(n)}\right] = \E \left[|\hat\mu_j - \mu_j|^{2\round[n]}\right] \le \binom{n}{\round[n]} \Sigma_{jj}^{\round[n]}.
\end{align*}
Plugging this into the bound \cref{eq:first_jensen_small_n} and using a similar approximation for the binomial term as before yields the final bound
\[
\E \left[\| \hat\mu - \mu \|_p \right] \le 2 \left( \sum_{j=1}^d\left(\Sigma_{jj}\right)^{\round[n]}\right)^{\frac{1}{2\round[n]}},
\]
as desired.
We now derive the uniform bound based on $s=\min\{ \bp ,\, 2\round[n]\}$.
Notice that in the case that $n \ge t(\bp)$ we have $s = \bp$ and we have the simple bound
\[
\E \left[\| \hat\mu - \mu \|_p \right] \le 2\frac{\left( \sum_{j=1}^d\left(\Sigma_{jj}\right)^{\frac{\bp}{2}}\right)^{\frac{1}{\bp}}}{\sqrt{\floor{n/t(\bp)}}} = 2\frac{\left( \sum_{j=1}^d\left(\Sigma_{jj}\right)^{\frac{s}{2}}\right)^{\frac{1}{s}}}{\sqrt{\floor{n/t(\bp)}}} \le 4\left( \sum_{j=1}^d\left(\Sigma_{jj}\right)^{\frac{s}{2}}\right)^{\frac{1}{s}}\sqrt{\frac{s}{n}},
\]
where the final inequality follows from the facts that $n \ge t(\bp)$ implies
$\floor{\frac{n}{t(\bp)}} \ge \frac{n}{2t(\bp)}$ and $
t(\bp)=2\left\lceil \frac{\bp}{2}\right\rceil-1 \le 2\bp = 2s$ and so $\frac{1}{\sqrt{\floor{n/t(\bp)}}} \le \sqrt{\frac{2t(\bp)}{n}} \le 2\sqrt{\frac{s}{n}}$.
Now in the case that $n < t(\bp)$ we have $s = 2\left\lceil \frac{n}{2}\right\rceil = r(n)$ and hence
\[
\E \left[\| \hat\mu - \mu \|_p \right] \le 2\left( \sum_{j=1}^d\left(\Sigma_{jj}\right)^{\round[n]}\right)^{\frac{1}{2\round[n]}} = 2\left( \sum_{j=1}^d\left(\Sigma_{jj}\right)^{\frac{s}{2}}\right)^{\frac{1}{s}} \le 2\left( \sum_{j=1}^d\left(\Sigma_{jj}\right)^{\frac{s}{2}}\right)^{\frac{1}{s}}\sqrt{\frac{s}{n}},
\]
where the final inequality uses the fact that $
s=2\left\lceil \frac{n}{2}\right\rceil \ge n$.
This completes the proof.
\end{proof}

\subsection{Lower bounds}

In this section we prove that the rate achieved in \cref{thm:medianofmeans} is optimal.
The lower bound has two components, proved as two separate lemmas, corresponding to the two regimes discussed in the introduction.
The first obstruction is already present in the Gaussian location model and gives the Gaussian-width branch.
The second obstruction is genuinely non-Gaussian and is witnessed by heavy-tailed distributions. We begin with the Gaussian lower bound. It matches \cref{thm:medianofmeans} up to constants in the regime
$\bp \le 2\round[n]$.

\begin{lemma}[Gaussian-width lower bound]\label{lem:lb_expectation_subgaussian_rate}
There exists a universal constant $c>0$ such that for any $d,n\in\mathbb N$, any $p\in[1,\infty]$, and any estimator
$\cA:(\mathbb R^d)^n\to\mathbb R^d$,
\[
\sup_{\mu\in\mathbb R^d}
\E_{X_1,\dots,X_n\sim \mc{N}(\mu,I_d)}
\bigl[
\|\cA(X_1,\dots,X_n)-\mu\|_{p}
\bigr]
\ge
c\,d^{\frac{1}{\bp}}\sqrt{\frac{\bp}{n}}.
\]
Consequently, the same lower bound holds for the minimax risk $\mathfrak R(n,d,p)$.
\end{lemma}

\begin{proof}
Let $P_\mu=\mc{N}(\mu,I_d)$ and let $G\sim \mc{N}(0,I_d)$.
Since the Gaussian family $\{P_\mu:\mu\in\mathbb R^d\}$ is contained in the class of distributions with covariance $I_d$, it suffices to prove the lower bound over this subfamily.
By \cite[Theorem~28.7]{PolyanskiyWu2025InformationTheory}, which is a consequence of Anderson's lemma \cite{Anderson1955SymmetricUnimodal}, applied with
$
\rho(x)=\|x\|_p,
$
the minimax risk in the Gaussian location model is
\[
\inf_{\widehat\mu}
\sup_{\mu\in\mathbb R^d}
\E_{X_1,\dots,X_n\sim P_\mu}
\bigl[
\|\widehat\mu(X_1,\dots,X_n)-\mu\|_p
\bigr]
=
\E\left\|\frac{G}{\sqrt n}\right\|_p
=
\frac{\E\|G\|_p}{\sqrt n}.
\]
Indeed, the function $\rho(x)=\|x\|_p$ is lower-semicontinuous and bowl-shaped (its sublevel sets are symmetric and convex).
Therefore every estimator $\cA$ satisfies
\begin{equation}\label{eq:expected_risk_norm_lb}
\sup_{\mu\in\mathbb R^d}
\E_{X_1,\dots,X_n\sim P_\mu}
\bigl[
\|\cA(X_1,\dots,X_n)-\mu\|_p
\bigr]
\ge
\frac{\E\|G\|_p}{\sqrt n}.
\end{equation}
It remains to lower bound the Gaussian width term. The standard Gaussian norm estimate gives
\[
\E\|G\|_p
\asymp
d^{\frac{1}{\bp}}\sqrt{\bp},
\qquad
\bp=\min\{ p, \, \max\{\log d, \, 1\} \}.
\]
Plugging this estimate into
\cref{eq:expected_risk_norm_lb} yields the claim.
\end{proof}

We now prove the second lower bound. Unlike \cref{lem:lb_expectation_subgaussian_rate}, this obstruction is not present in the Gaussian location model. It is a genuinely heavy-tailed obstruction and gives the small-sample branch of the minimax rate. Since $\|x\|_p\ge \|x\|_\infty$ for all $x\in\mathbb R^d$ and all $p\ge 1$, it suffices to prove a lower bound in the $\ell_\infty$ norm.

\begin{lemma}[Heavy-tailed small-sample lower bound]\label{lem:lowerboundexpectation_smalln}
There exists a universal constant $c>0$ satisfying the following.
For any $d,n\in\mathbb N$ with $\log d\ge 2\round[n]$ and any estimator
$\cA:(\mathbb R^d)^n\to\mathbb R^d$, there exists a distribution $P$ with mean
$\mu\in\mathbb R^d$ and covariance $I_d$ such that
\[
\E_{X_1,\dots,X_n\sim P}
\left[
\|\cA(X_1,\dots,X_n)-\mu\|_\infty
\right]
\ge
c\,d^{\frac{1}{2\round[n]}}.
\]
Consequently, in the regime $\log d\ge 2\round[n]$, the same lower bound holds for the minimax risk $\mathfrak R(n,d,p)$ for every $p\in[1,\infty]$.
\end{lemma}

\begin{proof}
We use the high-probability $\ell_\infty$ lower bound proven in
\cref{thm:lowerboundhighprob}.
Apply \cref{thm:lowerboundhighprob} with $\sigma=1$ and $\delta_0=1/4$.
Then for any estimator $\cA$, there exists a distribution $P$ on $\mathbb R^d$
with mean $\mu$ and covariance $I_d$ such that
\[
\Pr_{X_1,\ldots,X_n\sim P}
\left[
\left\|
\cA(X_1,\ldots,X_n)-\mu
\right\|_\infty
\ge
\frac{1}{4}
\left(
1-2d^{-\frac{1}{\round[n]}}
\right)
(4d)^{\frac{1}{2\round[n]}}
\right]
>
\frac14 .
\]
Here we used that, for $\delta_0=1/4$, the auxiliary parameter appearing in
\cref{thm:lowerboundhighprob} is
$
\left(\frac{4\delta_0}{d}\right)^{\frac{1}{\round[n]}}
=
d^{-\frac{1}{\round[n]}}.
$
By the assumption $\log d\ge 2\round[n]$, we have
$
d^{-\frac{1}{\round[n]}}\le e^{-2}.
$
Therefore,
\[
\frac{1}{4}
\left(
1-2d^{-\frac{1}{\round[n]}}
\right)
(4d)^{\frac{1}{2\round[n]}}
\ge
\frac{1-2e^{-2}}{4}
d^{\frac{1}{2\round[n]}}.
\]
Using the elementary inequality $\E Z\ge a\,\Pr[Z\ge a]$ for nonnegative random variables $Z$, we conclude that
\[
\E_{X_1,\dots,X_n\sim P}
\left[
\left\|
\cA(X_1,\ldots,X_n)-\mu
\right\|_\infty
\right]
\ge
\frac{1-2e^{-2}}{16}
d^{\frac{1}{2\round[n]}}.
\]
This proves the lemma with
$
c=\frac{1-2e^{-2}}{16}.
$
\end{proof}

\subsection{Proof of \cref{thm:expected_minimax_risk}}
We begin by showing that guarantee of \cref{alg:MOM} provided in \cref{thm:medianofmeans} achieves the claimed minimax risk in the isotropic case.
First, notice that in the isotropic case $\Sigma = I_d$, so the guarantee reduces to
\[
\mathfrak{R}(n,d,p)
\le
\begin{cases}
\displaystyle
2d^{\frac{1}{2\round[n]}},
&\bp > 2\round[n],\\[2ex]
\displaystyle
2\,\frac{d^{\frac{1}{\bp}}}{\sqrt{\left\lfloor n/t(\bp)\right\rfloor}},
& \bp \le 2\round[n].
\end{cases}
\]
Notice that guarantee of \cref{thm:medianofmeans} based on $s=\min\{ \bp ,\, 2\round[n]\}$ implies the upper bound
\[
\riskp \le 4d^{\frac{1}{s}}\sqrt{\frac{s}{n}}.
\]
We now combine the lower bounds from \cref{lem:lb_expectation_subgaussian_rate,lem:lowerboundexpectation_smalln}.
Consider the regime where $\bp > 2 \round[n]$. In this case, we certainly have that $\log d \ge 2\round[n]$, so \cref{lem:lowerboundexpectation_smalln} implies that
\[
\riskp \gtrsim d^{\frac{1}{2\round[n]}}.
\]
We also have the lower bound on the risk implied by \cref{lem:lb_expectation_subgaussian_rate} that tells us that
\[
\riskp \gtrsim d^{\frac{1}{\bp}}\sqrt{\frac{\bp}{n}}.
\]
Thus, all together we conclude the lower bound
\[
\riskp \gtrsim
\begin{cases}
 d^{\frac{1}{2\lceil n/2\rceil}}, &  \bp > 2\round[n],\\[1ex]
d^{\frac{1}{\bp}}\sqrt{\frac{\bp}{n}}, &  \bp \le 2\round[n],
\end{cases}
\]
which is equivalent to the following lower bound based on $s$:
\[
\riskp \gtrsim d^{\frac{1}{s}}\sqrt{\frac{s}{n}}.
\]
This concludes the proof. \qed

\section{Optimal High-probability Guarantees for $\ell_\infty$}
\label{sec:optimal_high_prob_ell_inf}

In this section, we formalize and prove \cref{thm:inf_high_prob_phase_transition} which characterizes the optimal high-probability estimation guarantees under the $\ell_\infty$ norm. First, define the following minimax risk:
\begin{gather*}
    r(n,d, \delta) = \inf_{\cA:(\mathbb R^d)^n\to\mathbb R^d} \lbrb{r > 0: \sup_{\substack{\mu \in \mb{R}^d \\ P:\, \E_P [X]=\mu\\ \operatorname{Cov}_P(X)=I_d}} \Pr_{X_1,\dots,X_n\sim P} \lsrs{\norm{\cA(X_1,\dots,X_n) - \mu}_\infty \geq r} \leq \delta} \numberthis{} \label{eq:high-prob-risk}.
\end{gather*}
The formalization of \cref{thm:inf_high_prob_phase_transition} is given below.
\begin{theorem}
    \label{thm:high_prob_phase_transition}
    There is a universal constant $c >0$ such that for any $n, d\in\mathbb{N}$, and $\delta \in (0, c)$, we have that:
    \[
    r(n, d, \delta) \asymp
    \begin{cases}
         \left(\dfrac d\delta\right)^{\frac{1}{2\round[n]}},
        &
        n \lesssim \log(d/\delta),
        \\[2ex]

        \sqrt{\dfrac{\log(d/\delta)}{n}},
        &
        n \gtrsim \log(d/\delta).
    \end{cases}
    \]
    Moreover, this rate is achieved by the coordinate-wise median-of-means estimator.
\end{theorem}
Note that the rates in the two cases coincide when $n = \Theta(\log (d / \delta))$. The proof is presented in two parts with the upper bound established in \cref{sec:upperboundinfinity} and the lower bound in \cref{sec:lowerboundinfinity}.

\subsection{High-probability Upper Bound in $\ell_\infty$}
\label{sec:upperboundinfinity}
In this section we prove high probability bounds for the performance of the coordinate wise median of means estimator for estimating the mean of a random vector in $ \mathbb{R}^{d} $ under the $ \ell_{\infty} $ norm.
Towards this end, we begin with the one dimensional case.
\begin{theorem}\label{thm:tail_upper_bound_scalar}
Let $X_1, \dots, X_n \in \mb{R}$ be i.i.d.\ random variables with finite mean $\mu$ and variance $\sigma^2$.
Fix a confidence parameter $\delta \in (0, 1)$.
Given the number of buckets $t = 2\lceil \log(1/\delta)\rceil$ and the sample $X_1, \dots, X_n $ as input, \cref{alg:MOM} returns an estimate $\hat\mu$ satisfying
\begin{equation}\label{eq:univariate_MoM_whp_upperbound}
|\hat\mu - \mu| \le \begin{cases}
2\sigma  \left(\frac{1}{\delta}\right)^{\frac{1}{2\round[n]}},& n < 2 \lceil \log\left(\frac{1}{\delta}\right)\rceil,
\\
4\sigma \sqrt{\frac{e}{n}\log \left(\frac{1}{\delta}\right)},& n \ge 2 \lceil \log\left(\frac{1}{\delta}\right)\rceil,
\end{cases}
\end{equation}
with probability at least $1-\delta$ over the randomness of $X_1, \dots, X_n$.
\end{theorem}

\begin{proof}

We split the analysis into two very different cases.
First, let us consider the case that the sample size $n$ is smaller than the number of buckets $t$.
In this case, \cref{alg:MOM} returns the median of the $n$ sample points.
Using Markov's inequality together with the moment bound of \cref{momentsmedian} with $q=2$ and $t = n$ implies that our output $\hat\mu$ satisfies
\[
\Pr\left[|\hat\mu - \mu| > 2\sigma \left(\frac{1}{\delta}\right)^{\frac{1}{2\round[n]}} \right] \le
\frac{ \delta\E\left[|\hat\mu - \mu|^{2\round[n]}\right]}{(2\sigma)^{2\round[n]}} \le \frac{\delta{n \choose \round[n]}}{ 2^{2\round[n]}}
< \frac{\delta 2^n}{2^{2\round[n]}}
\le \delta,
\]
where the second-to-last inequality uses a simple upper bound on centered binomial coefficients.

Now consider the (well studied) case that $n$ is at least $t$.
In this case, \cref{alg:MOM} returns the median of the $t$ bucket means.
Each bucket contains at least $\lfloor \frac{n}{t} \rfloor$ samples, so the variance of each bucket mean is at most $\frac{\sigma^2}{\lfloor n/t \rfloor}$.
Thus, Markov's inequality applied together with moment bound of \cref{momentsmedian} tells us that our output $\hat\mu = \mathrm{Median}(\hat\mu^{(1)}, \dots, \hat\mu^{(t)})$ satisfies
\[
\Pr\left[|\hat\mu - \mu| > 2\sigma \sqrt{e\left(\frac{1}{\lfloor n/t \rfloor}\right) } \right]
\le \left(\frac{\lfloor n/t \rfloor}{4e\sigma^2}\right)^{\round[t]} {t \choose \round[t]}\max_{i \in [t]}\E \left[(\hat\mu^{(i)} - \mu)^2\right]^{\round[t]} \le e^{-\round[t]} \le \delta.
\]

The above implies the second bound in \cref{eq:univariate_MoM_whp_upperbound} since $\floor{\frac{n}{t}} \ge \frac{n}{4\log(1/\delta)}$. This completes the proof.
\end{proof}
Applying the above median-of-means estimator to each coordinate with failure parameter $\delta/d$ together with a union bound gives us the following corollary for estimation under the $\ell_\infty$ norm.
\begin{corollary}\label{cor:tail_upper_bound_linfty}
Let $X_1, \dots, X_n \in \mb{R}^d$ be i.i.d.\ random vectors with finite mean $\mu$ and covariance $\sigma^2 I_d$.
Fix a confidence parameter $\delta \in (0, 1)$.
Given the number of buckets $t = 2\lceil \log(d/\delta)\rceil$ and the sample $X_1, \dots, X_n $ as input, \cref{alg:MOM} returns an estimate $\hat\mu$ satisfying
\begin{equation}
\|\hat\mu - \mu\|_{\infty}\le \begin{cases}
2\sigma  \left(\frac{d}{\delta}\right)^{\frac{1}{2\round[n]}},& n < 2 \lceil \log\left(\frac{d}{\delta}\right)\rceil,
\\
4\sigma \sqrt{\frac{e}{n}\log \left(\frac{d}{\delta}\right)},& n \ge 2 \lceil \log\left(\frac{d}{\delta}\right)\rceil,
\end{cases}
\end{equation}
with probability at least $1-\delta$ over the randomness of $X_1, \dots, X_n$.
\end{corollary}

\subsection{High-probability Lower Bound in $ \ell_{\infty} $}
\label{sec:lowerboundinfinity}
In this sub-section we show a high-probability lower bound for any estimator in $ \ell_{\infty} $. We prove our lower bound in two regimes: the small sample regime ($n \leq \log (d / \delta)$) in \cref{thm:lowerboundhighprob} and the large sample regime ($n \geq \log (d / \delta)$) in \cref{lem:gaussian_highprob_large_sample_lb}.

\begin{theorem}\label{thm:lowerboundhighprob}
  For any $\delta\in(0,1/2)$, any $ \sigma>0 $, any  $n,d\in\mathbb{N}$ and any deterministic estimator $\cA:(\mathbb{R}^d)^n\to\mathbb{R}^d$, there exists a distribution $P$ on $\mathbb{R}^d$ with mean $\mu$ and covariance matrix $\sigma^2 I_d$ such that
  \begin{align*}
   \p_{X_1, \dots, X_n \sim P}\left[\|\cA(X_1, \dots, X_n)-\mu\|_\infty \ge \frac{\sigma(1-2p)}{4}\left(\frac{d}{\delta}\right)^{\frac{1}{2\round[n]}}\right] > \delta,
  \end{align*}
  where $ p=(\frac{4\delta}{d})^{\frac{1}{\round[n]}} $.
\end{theorem}

\begin{proof}
If $p\ge 1/2$, then $1-2p\le 0$, so the lower bound in the theorem is
non-positive. In this case the claim is trivial since the norm is always non-negative and the claim follows for any distribution $P$ with covariance $\sigma^2 I_d$. Thus, we assume that $p<1/2$ throughout.
Define
\[
    \ell=\frac{\sigma}{\sqrt{p(1-p)}}.
\]
We consider two distributions $P_0$ and $P_1$ on $\{0,\ell\}$ given by
\[
    \Pr_{Y\sim P_0}[Y=\ell]=p,
    \qquad
    \Pr_{Y\sim P_1}[Y=\ell]=1-p.
\]
Equivalently, under $P_0$ the point $\ell$ occurs with probability
$p$, while under $P_1$ it occurs with probability $1-p$.
Write $\mu_0 =\mathbb E_{Y \sim P_0}[Y]=p\ell$ and  $\mu_1 =\mathbb E_{Y \sim P_1}[Y]=(1-p)\ell$ for the mean of $P_0$ and $P_1$ respectively.
It is immediate that both distributions have variance $p(1-p)\ell^2=\sigma^2$.
Moreover, we define the separation of the means to be
\[
    \Delta = \mu_1-\mu_0=(1-2p)\ell,
\]
and since $p<1/2$, we have $\Delta>0$.
We now construct the hard distribution. For each
$b=(b_1,\ldots,b_d)\in\{0,1\}^d$, let $P_b$ be the product distribution
on $\mathbb R^d$ whose $j$th coordinate has distribution $P_{b_j}$,
independently over $j=1,\ldots,d$.
Notice that $P_b$ has mean $\mu_b$ and covariance matrix $\sigma^2 I_d$ where $\mu_b = (\mu_{b_1},\ldots,\mu_{b_d})$.
We will show that for any deterministic estimator $\mc{A}$, there exists \emph{some} choice of $b \in \{0, 1\}^d$ witnessing the lower bound.
Given samples $X_1,\ldots,X_n\in\mathbb R^d$, define a decoder
$\widehat b=\widehat b(X_1,\ldots,X_n)\in\{0,1\}^d$ coordinate-wise by
\[
    \widehat b_j
    \in
    \arg\min_{u\in\{0,1\}}
    \left|
        \mc{A}(X_1,\ldots,X_n)_j-\mu_u
    \right|,
\]
with ties broken arbitrarily. Now, if $\left\|\mc{A}(X_1,\ldots,X_n)-\mu_b\right\|_\infty< \frac{\Delta}{2}$,
then every coordinate of $\mc{A}(X_1,\ldots,X_n)$ is strictly closer
to the correct mean than to the incorrect mean, hence
\begin{equation}\label{eq:implication_lb_small_sample}
    \left\|
        \mc{A}(X_1,\ldots,X_n)-\mu_b
    \right\|_\infty
    <
    \frac{\Delta}{2}
    \quad\Longrightarrow\quad
    \widehat b=b.
\end{equation}

Let $U$ be the uniform distribution over $\{0, 1\}^d$ and let $B=(B_1,\ldots,B_d) \sim U$.
Conditioned on $B=b$,
let $X_1,\ldots,X_n$ be i.i.d.\ samples from $P_b$. For each coordinate $j$, define
the coordinate samples
\[
    Z_j=(X_{1,j},X_{2,j},\ldots,X_{n,j})\in \mc{Z},
\]
where $\mc{Z} \coloneqq \{0,\ell\}^n$.
Conditioned on $B=b$, the random vectors $Z_1,\ldots,Z_d$ are independent
and $Z_j$ follows $P_{b_j}^{n}$. Finally define $
    \eta
    =
    \frac12
    \sum_{z\in\{0,\ell\}^n}
    \min\{P^{n}_0(z),P^{n}_1(z)\}$, which, by restricting the sum to those $z$ satisfying
$|\{i:z_i=\ell\}|=\floor{n/2}$, can be lower bounded as follows:
\begin{equation}\label{eq:eta_lowerbound}
    \eta\geq \frac12\binom{n}{\floor{\frac{n}{2}}}\min\{p^{\round[n]}(1-p)^{\floor{\frac{n}{2}}} ,(1-p)^{\round[n]}p^{\floor{\frac{n}{2}}} \}  \ge \frac12\binom{n}{\floor{\frac{n}{2}}}p^{\round[n]}2^{-\floor{\frac{n}{2}}} \ge \frac12 p^{\round[n]} = \frac{2\delta}{d},
\end{equation}
where the second inequality uses the fact that $p < 1/2$ twice, first to say the first term in the minimum is no larger than the second and then lower bound $ 1-p $ by $ 1/2 $, the final inequality uses the well known inequality $\binom nr \ge (\frac{n}{r})^r$ for $n \ge r \ge 1$ to conclude that  $2^{-\floor{\frac{n}{2}}}\binom{n}{\floor{\frac{n}{2}}} \ge 2^{-\floor{\frac{n}{2}}}(\frac{n}{\floor{\frac{n}{2}}})^{\floor{\frac{n}{2}}} \ge 1$,\footnote{If $n=1$, then the claimed lower bound trivially holds since $2^{-\floor{\frac{n}{2}}}= 1$ and $\binom{n}{\floor{\frac{n}{2}}} = \binom{1}{0} =1$ by convention.} and the equality follows from the definition of $p$.
We now claim that \emph{every} decoder
$\phi:\bigl(\{0,\ell\}^n\bigr)^d\to\{0,1\}^d$ satisfies
\[
    \Pr_{B \sim U, \; X_1, \dots, X_n \sim P_B}\left[\phi(Z)=B\right]
    \le
    (1-\eta)^d.
\]
Indeed, writing $Z = (Z_1, \dots, Z_d)$ and summing over the sample space we have that
\begin{align*}
     \Pr_{B \sim U, \; X_1, \dots, X_n \sim P_B}\left[\phi(Z)=B\right]
    &=
    2^{-d}
    \sum_{b\in\{0,1\}^d}
    \sum_{z \in \mc{Z}^d}
    \mathbf 1\{\phi(z)=b\}
    \prod_{j=1}^d P^{n}_{b_j}(z_j) \\
    &=
    2^{-d}
    \sum_{z \in \mc{Z}^d}
    \prod_{j=1}^d P^{n}_{\phi(z)_j}(z_j) \\
    &\le
    2^{-d}
    \sum_{z \in \mc{Z}^d}
    \prod_{j=1}^d \max_{u\in\{0,1\}} P^{n}_{u}(z_j) \\
    &=
    2^{-d}
    \sum_{z_1\in\mc{Z}}\cdots\sum_{z_d\in\mc{Z}}
    \prod_{j=1}^d \max_{u\in\{0,1\}} P^{n}_{u}(z_j) \\
    &=
    \left(
        \frac12
        \sum_{w\in\mc{Z}}
        \max_{u\in\{0,1\}} P^{n}_{u}(w)
    \right)^d.
\end{align*}
Now notice that
\[
    \frac12
    \sum_{w \in \mc{Z}}
    \max\{P^{n}_{0}(w),P^{n}_{1}(w)\}
    =
    \frac12
    \sum_{w \in \mc{Z}}\left( P^{n}_{0}(w) + P^{n}_{1}(w) -
    \min\{P^{n}_{0}(w),P^{n}_{1}(w)\}\right)
    =
    1-\eta,
\]
so,
\[
    \Pr_{B \sim U, \; X_1, \dots, X_n \sim P_B}\left[\phi(Z)=B\right]\le (1-\eta)^d.
\]
The above upper bound certainly holds for our decoder $\widehat b$, so an averaging argument implies the existence of \emph{some} $b^\star \in \{0, 1\}^d$ such that
\[
    \Pr_{X_1, \dots, X_n \sim P_{b^\star}}\left[\widehat b(Z)=b^\star \right]\le (1-\eta)^d \le \left(1-\frac{2\delta}{d}\right)^d \le e^{-2\delta} < 1-\delta,
\]
where we used the lower bound on $\eta$ from~\cref{eq:eta_lowerbound} and the inequality $e^{-2\delta} \le 1-\delta$ that holds for $\delta \in (0,\frac{1}{2})$.
Now, using the above with \cref{eq:implication_lb_small_sample}, we conclude that
\[
    \Pr_{X_1, \dots, X_n \sim P_{b^\star}}
    \left[
        \left\|
            \mc{A}(X_1,\ldots,X_n)-\mu_{b^\star}
        \right\|_\infty
        \ge
        \frac{\Delta}{2}
    \right]
     \ge \Pr_{X_1, \dots, X_n \sim P_{b^\star}}\left[\widehat b(Z)\not=b^\star\right] > \delta.
\]
All that remains is to lower bound $\Delta/2$.
We have
\begin{align*}
    \frac{\Delta}{2}
    &=
    \frac{\sigma(1-2p)}{2\sqrt{p(1-p)}} \\
    &\ge
    \frac{\sigma(1-2p)}{2}
    \left(\frac d{4\delta}\right)^{\frac{1}{2\round[n]}} \\
    &\ge
    \frac{\sigma(1-2p)}{4}
    \left(\frac d\delta\right)^{\frac{1}{2\round[n]}}.
\end{align*}
where we used the facts $ p= \left(\frac{4\delta}{d}\right)^{\frac{1}{\round[n]}} $,  $\sqrt\frac{1}{1-p} \ge 1$ and $4^{-\frac{1}{2\round[n]}} \ge \frac{1}{2}$.
This completes the proof.
\end{proof}
\begin{lemma}[Gaussian lower bound for the large-sample regime]
\label{lem:gaussian_highprob_large_sample_lb}
There exists a universal constant \(c>0\) such that for every
\(n,d\in\mathbb N\), every \(\sigma>0\), and every
\(\delta\in(0,1/4)\),
\[
    r(n,d,\sigma, \delta)
    \ge
    c\sigma\sqrt{\frac{\log(d/\delta)}{n}} .
\]
\end{lemma}

\begin{proof}
Note that it suffices to prove the lower bound on the Gaussian distribution:
\[
    X_1,\ldots,X_n \mid \mu
    \sim
    \mathcal N(\mu,\sigma^2 I_d)^{\otimes n},
    \qquad
    \mu\in\mathbb R^d.
\]

Fix an arbitrary estimator
\[
    \cA:(\mathbb R^d)^n\to\mathbb R^d.
\]
Let
\[
    \bar X=\frac1n\sum_{i=1}^n X_i.
\]
For a parameter \(s>0\) and assume  the Gaussian prior on $\mu$:
\[
    \mu\sim \mathcal N(0,s^2 I_d).
\]
Under this prior and the Gaussian observation model above, \(\bar X\) is a
sufficient statistic for \(\mu\). Concretely, the posterior distribution
of \(\mu\) given the full sample \(\bm{X} =(X_1,\ldots,X_n)\) is
\[
    \mu \mid \bm{X}
    \sim
    \mathcal N\left(
        \alpha_s \bar X,\,
        \tau_s^2 I_d
    \right),
\]
where
\[
    \alpha_s
    =
    \frac{s^2}{s^2+\sigma^2/n},
    \qquad
    \tau_s^2
    =
    \left(\frac1{s^2}+\frac n{\sigma^2}\right)^{-1}
    =
    \frac{s^2\sigma^2}{\sigma^2+n s^2}.
\]
Note, in particular, that:
\[
    \tau_s \longrightarrow \frac{\sigma}{\sqrt n}
    \qquad
    \text{as } s\to\infty.
\]
Fix \(r>0\), and define the cube
\[
    K_r
    \coloneqq
    \{z\in\mathbb R^d:\|z\|_\infty<r\}.
\]
Observe that $K_r$ is convex and centrally symmetric.

Now condition on the realized sample \(\bm{X}\). The conditional success probability of the estimator is
\[
    \Pr\left[
        \|\cA(\bm{X})-\mu\|_\infty<r
        \,\middle|\,
        \bm{X}
    \right]
    =
    \Pr\left[
        \mu\in \cA(\bm{X}) + K_r
        \,\middle|\,
        \bm{X}
    \right].
\]
Since
\[
    \mu\mid \bm{X}
    \sim
    \mathcal N(\alpha_s \bar{X},\tau_s^2 I_d),
\]
this conditional success probability of the estimator is
\[
    \Pr\left[
        \alpha_s \bar{X} + \tau_s G \in \cA (\bm{X}) + K_r \,\middle|\,
        \bm{X}
    \right]
    =
    \Pr\left[
        \tau_s G \in \cA (\bm{X}) - \alpha_s \bar{X} + K_r \,\middle|\,
        \bm{X}
    \right],
\]
where \(G\sim\mathcal N(0,I_d)\). By Anderson's inequality, for a centered Gaussian vector and a convex centrally symmetric set \(K_r\), the Gaussian measure of a translate of
\(K_r\) is maximized at the centered translate \cite{Anderson1955SymmetricUnimodal}. Hence
\[
    \Pr\left[
        \tau_s G \in \cA (\bm{X}) - \alpha_s \bar{X} + K_r \,\middle|\,
        \bm{X}
    \right]
    \le
    \Pr\left[
        \tau_s G\in K_r \,\middle|\,
        \bm{X}
    \right].
\]
Averaging over \(X\) under the joint prior-predictive distribution gives
\[
    \Pr\left[
        \|\cA(\bm X)-\mu\|_\infty<r
    \right]
    \le
    \Pr\left[
        \tau_s \|G\|_\infty<r
    \right].
\]
And consequently,
\[
    \Pr\left[
        \|\cA(\bm X)-\mu\|_\infty\ge r
    \right]
    \ge
    \Pr\left[
        \tau_s \|G\|_\infty\ge r
    \right].
\]
Letting \(s\to\infty\), and using continuity of the distribution of
\(\|G\|_\infty\), yields
\begin{equation*}
    \sup_{\mu\in\mathbb R^d}
    \Pr_{X_1,\ldots,X_n\sim \mathcal N(\mu,\sigma^2 I_d)}
    \left[
        \|\cA(X_1,\ldots,X_n)-\mu\|_\infty\ge r
    \right]
    \ge
    \Pr\left[
        \frac{\sigma}{\sqrt n}\|G\|_\infty\ge r
    \right]. \numberthis{} \label{eq:lb_large_sample}
\end{equation*}
It remains to choose \(r\). We use the following standard Gaussian maximum
lower bound:
\begin{claim}
    \label{clm:gau_max_lb}
    There exists a universal constant \(c_0>0\) such that for every \(d\in\mathbb N\) and every \(\delta\in(0,1/4)\),
\[
    \Pr\left[
        \|G\|_\infty
        \ge
        c_0\sqrt{\log(d/\delta)}
    \right]
    \ge
    \delta .
\]
\end{claim}
\begin{proof}
    Let
    \[
        L=\log(d/\delta),
        \qquad
        t=c_0\sqrt L,
        \qquad
        p_t=\Pr(|g|\ge t),
    \]
    where \(g\sim\mathcal N(0,1)\). If \(t\le 1\), then
    \[
        \Pr[\|G\|_\infty\ge t]
        \ge
        \Pr[|g|\ge 1]
        >
        \delta,
    \]
    as $ \Pr[|g|\ge 1]\geq 1/4 $ and \(\delta<1/4\).

    If \(t>1\), the Gaussian tail lower bound gives
    \[
        p_t
        \ge
        c_1 t^{-1}e^{-t^2/2}
    \]
    for a universal constant \(c_1>0\). Thus
    \[
        d p_t
        \ge
        \frac{c_1 d}{c_0\sqrt L}
        \exp\left(-\frac{c_0^2}{2}L\right)
        =
        \frac{c_1\delta}{c_0\sqrt L}
        \left(\frac d\delta\right)^{1-c_0^2/2}.
    \]
    Choosing \(c_0>0\) sufficiently small ensures that \(d p_t\ge 2\delta\) for
    all \(d\ge 1\) and \(\delta\in(0,1/4)\).
     Using the independence of the coordinates yields,
    \[
        \Pr[\|G\|_\infty\ge t]
        =
        1-(1-p_t)^d
        \ge
        1-e^{-d p_t}
        \ge
        1-e^{-2\delta}
        \ge
        \delta .
    \]
    This proves the Gaussian maximum lower bound.
\end{proof}

Now setting
\[
    r
    =
    c_0\sigma\sqrt{\frac{\log(d/\delta)}{n}}.
\]
with \cref{clm:gau_max_lb} and Equation~\ref{eq:lb_large_sample} concludes the proof.
\end{proof}

\subsection*{Proof of \cref{thm:high_prob_phase_transition}}

Write $L=\log(d/\delta)$.
We first prove the upper bound. In the
isotropic case, \cref{cor:tail_upper_bound_linfty} gives, for the
coordinate-wise median-of-means estimator with $t=2\ceil{L}$ buckets
\[
r(n,d,\delta)
\lesssim
\begin{cases}
\left(\dfrac d\delta\right)^{\frac{1}{2\round[n]}},
& n<2\ceil{L},
\\[2ex]
\sqrt{\dfrac{L}{n}},
& n\ge 2\ceil{L}.
\end{cases}
\]
This implies the claimed upper bound in both regimes. If $n\lesssim L$, then
the first line is exactly the desired bound when it applies, while if the
second line applies then $\sqrt{L/n}\le 1\le
(d/\delta)^{1/(2\round[n])}$. If $n\gtrsim L$, then the second line is exactly
the desired bound when it applies, while if the first line applies then
$n=\Theta(L)$, and the two rates are comparable up to universal constants.

We now prove the lower bound. First, consider the large sample regime $n \gtrsim L$.
We can invoke the Gaussian lower bound~(\cref{lem:gaussian_highprob_large_sample_lb})
in the isotropic to get
\[
r(n,d,\delta)
\gtrsim
\sqrt{\frac{L}{n}},
\]
as claimed.

It remains to prove the lower bound in the small sample regime.
Assume that $n\lesssim L$.
By taking the universal constant small in the theorem statement to be small enough, we may assume that $\delta<1/4$ and that $L\ge 2\log 4$. Set
\[
a=\left(\frac{4\delta}{d}\right)^{1/\round[n]}.
\]
If $a\le 1/4$, then $1-2a\ge 1/2$, and
\cref{thm:lowerboundhighprob} implies that for every estimator there is a
distribution with mean $\mu$ and covariance $I_d$ such that
\[
\Pr\left[
\|\cA(X_1,\ldots,X_n)-\mu\|_\infty
\ge
\frac18
\left(\frac d\delta\right)^{\frac{1}{2\round[n]}}
\right]
>
\delta.
\]
Hence, by the definition of the minimax risk, we have
\[
r(n,d,\delta)
\gtrsim
\left(\frac d\delta\right)^{\frac{1}{2\round[n]}}.
\]

On the other hand if $a>1/4$, then
\[
\left(\frac d\delta\right)^{\frac{1}{2\round[n]}}
\asymp 1,
\]
which we can combine with \cref{lem:gaussian_highprob_large_sample_lb} and our assumption that $n \lesssim L$ to conclude that
\[
r(n,d,\delta)
\gtrsim \sqrt{\frac{L}{n}} \gtrsim  1 \asymp
\left(\frac d\delta\right)^{\frac{1}{2\round[n]}}.
\]
Combining the upper and lower bounds completes the
proof.
\qed

\paragraph{Funding.}
Mikael was while carrying out this work, supported by an Internationalisation Fellowship from the Carlsberg Foundation, Independent Research Fund Denmark (DFF) Sapere Aude Research
Leader grant No. 9064-00068B, and the European Union (ERC, TUCLA, 101125203).
\newline
Kasper was supported by the European Union (ERC, TUCLA, 101125203). Views and opinions expressed are however
those of the author(s) only and do not necessarily reflect those of the European Union or the European
Research Council. Neither the European Union nor the granting authority can be held responsible for them.
\newline
 Yeshwanth Cherapanamjeri is supported by an NTU-NAP Start-up Grant (025661-00016).

\paragraph{AI assistance disclosure.}
The theorem statements, proofs, and hard constructions in this work were developed by the authors. Most of the results were obtained in summer 2025 without material AI assistance. During the preparation of the manuscript in April 2026, OpenAI GPT-5.5 assisted with formalizing and drafting parts of the proof of the lower bound for projection based median of means estimators, as well as with checking and polishing the manuscript.

\bibliographystyle{alpha}
\bibliography{refs.bib}
\appendix

\section{Omitted Proofs}\label{app:omitted}

\begin{lemma}\label{lem:integralproduct}
    Let $ f_1, \ldots, f_t:\mathbb{R}_{+}\rightarrow\mathbb{R}_{+} $ be $t$, measurable, non-increasing and integrable functions. Then
    \begin{align*}
    \int_{0}^{\infty} \prod_{i=1}^{t} f_{i}(x^{1/t})  dx \leq \prod_{i=1}^{t} \int_{0}^{\infty} f_{i}(x) dx
    \end{align*}
\end{lemma}

\begin{proof}
  We have that
  \begin{align*}
   \prod_{i=1}^{t} \int_{0}^{\infty} f_{i}(x) dx
   &=
    \prod_{i=1}^{t} \int_{0}^{\infty} f_{i}(x_{i}) dx_{i}
   \\
   &=\int_{0}^{\infty}\ldots \int_{0}^{\infty} \prod_{i=1}^{t} f_{i}(x_{i}) dx_{1}\ldots dx_{t}
   \\
   &= \sum_{\stackrel{x_{i_{1}}<\ldots <x_{i_{t}}}{\left\{  i_{1},\ldots ,i_{t}\right\} }=[t]} \int_{0}^{\infty}\int_{0}^{x_{i_{t}}}\ldots \int_{0}^{x_{i_{2}}} \prod_{i=1}^{t} f_{i}(x_{i}) dx_{i_{1}}\ldots dx_{i_{t-1}} dx_{i_{t}}
   \\
   &\geq
   \sum_{\stackrel{x_{i_{1}}<\ldots <x_{i_{t}}}{\left\{  i_{1},\ldots ,i_{t}\right\} }=[t]} \int_{0}^{\infty} \prod_{i=1}^{t} f_{i}(x_{i_{t}}) \int_{0}^{x_{i_{t}}} \ldots \int_{0}^{x_{i_{2}}} 1 dx_{i_{1}}\ldots dx_{i_{t-1}} dx_{i_{t}}
   \\
   &=
   t! \int_{0}^{\infty} \prod_{i=1}^{t} f_{i}(x) \frac{x^{t-1}}{(t-1)!} dx
   \\
   &=\int_{0}^{\infty} \prod_{i=1}^{t} f_{i}(x^{1/t}) dx,
  \end{align*}
  where the third equality follows from decomposing the integration domain $ \mathbb{R}_{+}^{t} $ into the $ t! $ disjoint regions corresponding to the possible orderings $ x_{i_{1}} < \cdots < x_{i_{t}} $ of the coordinates. The inequality follows from the fact that each $ f_i $ is non-increasing, so on any region where $ x_{i_{1}} < \cdots < x_{i_{t}} $ we have $ f_i(x_i) \geq f_i(x_{i_t}) $ for every $ i \in [t] $. The second last equality follows from the fact that there are $ t! $ such orderings, together with an explicit evaluation of the inner integrals. The final equality follows from a change of variables.
\end{proof}

\momentsmedian*

\begin{proof}
  The case $ q=0 $ follows from the left hand side of the inequality to be less than or equal to the right hand side of the inequality. Furthermore, if the right hand side is infinite then the inequality also holds. Thus, we can assume that $ q>0 $ and that the right hand side is finite, meaning all the expectations $ \e\left[ |X_{i}-\mu|^{q}\right] $ are finite.
  By calculations we have that
  \begin{align*}
    \e\left[|M-\mu|^{ q \round[t]}\right]
    &=\int_{0}^{\infty} \p\left[|M-\mu|^{ q \round[t]}\geq x\right] dx
    =\int_{0}^{\infty} \p\left[|M-\mu|^{q}\geq x^{\frac{1}{\round[t]}}\right] dx,
  \end{align*}
  where we have used that the expectation of a positive random variable is the integral of its right tail.
  Now for the median $ M $ of $ X_{1},\ldots,X_{t} $ to be $ x $ away from the mean $ \mu $ it must be the case that there are $ \round[t] $ $X_i$'s such that $ |X_{i}-\mu|\geq x.$ By this, a union bound and independence of $ X_{i} $'s we then get that:
  \begin{align*}
    \e\left[|M-\mu|^{ q \round[t]}\right]
    &\leq \int_{0}^{\infty} \p\left[\exists I \in \binom{t}{\round[t]}: \forall i \in I \text{ such that } |X_{i}-\mu|^{q}\geq x^{\frac{1}{\round[t]}}\right] dx
    \\
    &\leq  \sum_{I\in \binom{t}{\round[t]}} \int_{0}^{\infty} \p\left[\forall i \in I: |X_{i}-\mu|^{q}\geq x^{\frac{1}{\round[t]}}\right] dx
    \\
    &=
     \sum_{I\in \binom{t}{\round[t]}} \int_{0}^{\infty} \prod_{i\in I}\p\left[ |X_{i}-\mu|^{q}\geq x^{\frac{1}{\round[t]}}\right] dx.
  \end{align*}
Now using that for each $ i $, the function $\p\left[ |X_{i}-\mu|^{q}\geq x\right]$ is non-increasing in $ x $, integrable by considering the case where the expectation is finite for all $ X_{i} $, and measureable by $ X_{i} $ being a random variable, we can use \cref{lem:integralproduct} to get that
\begin{align*}
\e\left[|M-\mu|^{ q \round[t]}\right]
&\leq
\sum_{I\in \binom{t}{\round[t]}} \prod_{i\in I} \int_{0}^{\infty} \p\left[ |X_{i}-\mu|^{q}\geq x\right] dx
\\
&=\sum_{I\in \binom{t}{\round[t]}} \prod_{i\in I} \e\left[ |X_{i}-\mu|^{q}\right]
\\
&\leq  \binom{t}{\round[t]} \max_{i\in[t]} \e\left[ |X_{i}-\mu|^{q}\right]^{\round[t]}.
\end{align*}
where the first inequality follows from \cref{lem:integralproduct}, the first equality follows from the fact that the expectation of a positive random variable is the integral of its right tail, and the last inequality follows from taking max under the product.
\end{proof}

\section{Proof of \cref{thm:proj_sub_optimal}}
\label{ssec:app_sub_opt_proj_median}

In this section, we establish that projection-based Median-of-Means estimators are provably sub-optimal for $\ell_p$ norms when $p > 2$. Let $p > 2$ and $q < 2$ denote the order of the corresponding dual ball. The projection-based median is now defined as follows for a pointset $Y = \{y_i\}_{i = 1}^n$ with $y_i \in \mathbb{R}^d$:

\begin{equation*}
    \medproj (Y) = \argmin_x \max_{v: \norm{v}_q = 1} \abs*{\inp{v}{x} - \med \lprp{\lbrb{\inp{v}{y_i}}_{i = 1}^n}}. \tag{PROJ-MED} \label{eq:proj_med}
\end{equation*}
\ProjSuboptimal*
\begin{remark}
    \label{rmk:proj_sub_optimal}
    Firstly, observe that the lower bound in \cref{thm:proj_sub_optimal} directly translates to an in-expectation lower bound for the projection-based median-of-means estimator. Secondly, observe that \cref{thm:medianofmeans} when $\Sigma = I_d$ (and $n = \omega(\log (d))$) yields the in-expecation bound:
    \begin{equation*}
        \E \lsrs{\norm{\wh{\mu} - \mu}_p} \leq 4 d^{1 / p} \sqrt{\frac{\log d}{n}}.
    \end{equation*}
    Comparing this bound with \cref{thm:proj_sub_optimal}, we see that the projection-based median-of-means estimator is strictly worse when $n = \widetilde{o} (\sqrt{d})$. In particular, for $p = \infty$, \cref{thm:proj_sub_optimal} incurs a \emph{polynomial} dependence on $d$ while \cref{thm:medianofmeans} is \emph{logarithmic}.
\end{remark}

\newcommand{\numbuckets}{t}
\newcommand{\bucksize}{m}

We will show that the median-of-means estimator with the projection-based median defined above is sub-optimal for any choice of number of buckets, $\numbuckets$, equivalently bucket size $\bucksize = n / \numbuckets$, for the distribution $D$, corresponding to the uniform distribution over the signed and scaled basis vectors. Formally:
\begin{equation*}
    \forall i \in [d], s \in {\pm 1}: \Pr_{X \thicksim D} \lprp{X = s \sqrt{d} e_i} = \frac{1}{2d}.
\end{equation*}
where $ e_{i} $ denote the standard basis vectors in $\mathbb{R}^d$. Note that $D$ has mean $0$ and covariance $I_d$.

Now, let $\bm{X} = \{X_i\}_{i = 1}^n$ be $n$ iid samples from from $D$. Furthermore, from our definition of $D$, $X_i = s_i \sqrt{d} e_{c(i)}$ with $s_i = \pm 1$ and $c(i) \in [d]$. We first prove a simple lemma establishing that the $c(i)$ are unique with high-probability:
\begin{lemma}
    \label{lem:unique_coords}
    Under the preconditions of \cref{thm:proj_sub_optimal}, we have:
    \begin{equation*}
        \forall i, j\in [n], i \neq j: c(i) \neq c(j)
    \end{equation*}
    with probability at least $0.9$.
\end{lemma}
\begin{proof}
    We have:
    \begin{align*}
        \p \lprp{\forall i, j\in [n], i \neq j: c(i) \neq c(j)} &= \prod_{i = 1}^n \lprp{1 - \frac{(i - 1)}{d}} \geq \prod_{i = 1}^n \lprp{1 - \frac{n}{d}} \\
        &= \lprp{1 - \frac{n}{d}}^n \geq \exp\lprp{- \frac{2n}{d} \cdot n} \geq 0.9
    \end{align*}
    as $n \leq c \sqrt{d}$.
\end{proof}
Now, we condition on the uniqueness of $c(i)$. Define the bucketed means $Z_i$ as follows for $i \in [t]$ with $m \coloneqq n / t$:
\begin{equation*}
    Z_i = \frac{1}{m} \sum_{j = (i - 1)m + 1}^{i m} X_j
\end{equation*}
with the collection denoted by $\bm{Z}$. From \ref{eq:proj_med}, the projection-based median-of-means estimate is the minimizer of the following function:
\begin{equation*}
    F(x) \coloneqq \max_{v: \norm{v}_q = 1} \abs*{\inp{v}{x} - \med \lprp{\lbrb{\inp{v}{Z_i}}_{i = 1}^t}}
\end{equation*}
We now make some observations:
\begin{observation}
    \label{obs:f_conv}
    $F$ is a convex function of $x$.
\end{observation}
\begin{proof}
    This follows from $F$ being defined as the maximum of a set of convex functions.
\end{proof}
Next, $F$ satisfies the following symmetries:
\begin{observation}
    \label{obs:sym_f}
    Defining:
    \begin{equation*}
        \wt{X}_i = \sqrt{d} e_{c (i)}, S \in \mathbb{R}^{d \times d}: S_{i, j} \coloneqq
        \begin{cases}
            0 & \text{if } i \neq j \\
            1 & \text{if } i=j, \forall \ell, c(\ell) \neq i \\
            1 & \text{if } i=j,  \exists \ell, c(\ell) = i \text{ and } s_\ell = 1 \\
            -1 & \text{otherwise}
        \end{cases}, \text{ and }
        \wt{Z}_i \coloneqq \frac{1}{m} \sum_{j = (i - 1)m + 1}^{i m} \wt{X}_j.
    \end{equation*}
    where uniqueness of $c(i)$ from \cref{lem:unique_coords} ensures that $S$ is well-defined.
    Finally, defining $\wt{F}$ analogously to $F$, we observe the following:
    \begin{equation*}
        \forall x: F(x) = \wt{F} (S x).
    \end{equation*}
    Consequently,
    \begin{equation*}
        \argmin_{x} F(x) = S \argmin_{x} \wt{F} (x).
    \end{equation*}
\end{observation}
\begin{proof}
    From the definition of $S$ we have $$ S\tilde{X}_{i}=S\sqrt{d}e_{c(i)}=S_{c(i),c(i)}\sqrt{d}e_{c(i)}=s_{i}\sqrt{d}e_{c(i)}=X_{i} $$ which implies
    \begin{equation*}
        \forall v, \norm{v}_q = 1, i \in [t]: \inp{v}{Z_i} = \inp{Sv}{\wt{Z}_i}.
    \end{equation*}
    and hence
    \begin{align*}
      F(x) = \max_{v: \norm{v}_q = 1} \abs*{\inp{v}{x} - \med \lprp{\lbrb{\inp{v}{Z_i}}_{i = 1}^t}}=  \max_{v: \norm{v}_q = 1} |\inp{Sv}{Sx} - \med \lprp{\lbrb{\inp{Sv}{\tilde{Z}_i}}_{i = 1}^t}|= \tilde{F}(Sx)
    \end{align*}
    where the last equality follows from the fact that \(S\) is a diagonal matrix with entries in \(\{\pm 1\}\), so \(Sv\) ranges over all vectors with \(\ell_q\) norm at most \(1\) as \(v\) does.
\end{proof}
On the event of \cref{lem:unique_coords}, define the bucket supports
\begin{equation*}
    C_r \coloneqq \{c((r - 1)m + 1), \dots, c(rm)\}, \qquad r \in [t].
\end{equation*}
Then the sets \(C_1, \dots, C_t\) are disjoint, \(|C_r| = m\) for all \(r \in [t]\), and
\begin{equation*}
    \wt{Z}_r = \frac{\sqrt{d}}{m} \sum_{j \in C_r} e_j.
\end{equation*}
Hence,
\begin{equation*}
    \wt{F}(x)
    =
    \max_{v: \norm{v}_q = 1}
    \abs*{
        \inp{v}{x}
        -
        \med \lprp{
            \lbrb{
                \frac{\sqrt{d}}{m} \sum_{j \in C_r} v_j
            }_{r = 1}^t
        }
    }.
\end{equation*}

\begin{lemma}
    \label{lem:within_bucket_symmetry}
    Fix \(r \in [t]\), and let \(P\) be any permutation matrix that permutes the coordinates in \(C_r\) and fixes all coordinates outside \(C_r\). Then
    \begin{equation*}
        \forall x \in \mathbb{R}^d: \wt{F}(Px) = \wt{F}(x).
    \end{equation*}
\end{lemma}
\begin{proof}
    Since
    \begin{equation*}
        \wt{Z}_\ell = \frac{\sqrt{d}}{m} \sum_{j \in C_\ell} e_j
    \end{equation*}
    for every \(\ell \in [t]\), permuting coordinates inside \(C_r\) leaves every \(\wt{Z}_\ell\) unchanged:
    \begin{equation*}
        \forall \ell \in [t]: P \wt{Z}_\ell = \wt{Z}_\ell.
    \end{equation*}
    Therefore,
    \begin{align*}
        \wt{F}(Px)
        &=
        \max_{v: \norm{v}_q = 1}
        \abs*{
            \inp{v}{Px}
            -
            \med \lprp{
                \lbrb{
                    \inp{v}{\wt{Z}_\ell}
                }_{\ell = 1}^t
            }
        } \\
        &=
        \max_{v: \norm{v}_q = 1}
        \abs*{
            \inp{P^\top v}{x}
            -
            \med \lprp{
                \lbrb{
                    \inp{P^\top v}{\wt{Z}_\ell}
                }_{\ell = 1}^t
            }
        }.
    \end{align*}
    Since \(P\) is a permutation matrix, \(\norm{P^\top v}_q = \norm{v}_q\). Hence, after the change of variables \(u = P^\top v\),
    \begin{equation*}
        \wt{F}(Px)
        =
        \max_{u: \norm{u}_q = 1}
        \abs*{
            \inp{u}{x}
            -
            \med \lprp{
                \lbrb{
                    \inp{u}{\wt{Z}_\ell}
                }_{\ell = 1}^t
            }
        }
        =
        \wt{F}(x).
    \end{equation*}
\end{proof}

\begin{corollary}
    \label{cor:bucketwise_constant}
    There exists a minimizer \(x^\star\) of \(\wt{F}\) such that for every \(r \in [t]\), all coordinates of \(x^\star\) in \(C_r\) are equal.
\end{corollary}
\begin{proof}
    Let \(\mathcal{G}_{\mathrm{in}}\) denote the finite group of permutation matrices that permute coordinates inside each \(C_r\) and fix all coordinates outside \(\cup_{r = 1}^t C_r\). For any minimizer \(x\) of \(\wt{F}\), define its symmetrization
    \begin{equation*}
        \bar{x} \coloneqq \frac{1}{|\mathcal{G}_{\mathrm{in}}|} \sum_{P \in \mathcal{G}_{\mathrm{in}}} P x.
    \end{equation*}
    By convexity of \(\wt{F}\) from \cref{obs:f_conv,lem:within_bucket_symmetry},
    \begin{equation*}
        \wt{F}(\bar{x})
        \leq
        \frac{1}{|\mathcal{G}_{\mathrm{in}}|} \sum_{P \in \mathcal{G}_{\mathrm{in}}} \wt{F}(Px)
        =
        \frac{1}{|\mathcal{G}_{\mathrm{in}}|} \sum_{P \in \mathcal{G}_{\mathrm{in}}} \wt{F}(x)
        =
        \wt{F}(x).
    \end{equation*}
    Since \(x\) is already a minimizer, \(\bar{x}\) is also a minimizer.

    By construction, \(\bar{x}\) is invariant under every permutation inside each \(C_r\). Hence,
    \begin{equation*}
        \forall r \in [t], \; \forall i, j \in C_r: \bar{x}_i = \bar{x}_j.
    \end{equation*}
\end{proof}
\begin{lemma}
    \label{lem:bucket_swap_symmetry}
    Let \(\sigma\) be any permutation of \([t]\). For each \(r \in [t]\), choose any bijection
    \begin{equation*}
        \phi_r : C_r \to C_{\sigma(r)}.
    \end{equation*}
    Let \(P_{\sigma,\phi}\) be the permutation matrix that sends each coordinate \(j \in C_r\) to \(\phi_r(j)\), and fixes all coordinates outside \(\cup_{r = 1}^t C_r\). Then
    \begin{equation*}
        \forall x \in \mathbb{R}^d: \wt{F}(P_{\sigma,\phi} x) = \wt{F}(x).
    \end{equation*}
\end{lemma}
\begin{proof}
    By construction,
    \begin{equation*}
        P_{\sigma,\phi} \wt{Z}_r
        =
        \frac{\sqrt{d}}{m} \sum_{j \in C_r} P_{\sigma,\phi} e_j
        =
        \frac{\sqrt{d}}{m} \sum_{k \in C_{\sigma(r)}} e_k
        =
        \wt{Z}_{\sigma(r)}.
    \end{equation*}
    Therefore,
    \begin{align*}
        \wt{F}(P_{\sigma,\phi} x)
        &=
        \max_{v: \norm{v}_q = 1}
        \abs*{
            \inp{v}{P_{\sigma,\phi} x}
            -
            \med \lprp{
                \lbrb{
                    \inp{v}{\wt{Z}_r}
                }_{r = 1}^t
            }
        } \\
        &=
        \max_{v: \norm{v}_q = 1}
        \abs*{
            \inp{P_{\sigma,\phi}^\top v}{x}
            -
            \med \lprp{
                \lbrb{
                    \inp{P_{\sigma,\phi}^\top v}{\wt{Z}_{\sigma(r)}}
                }_{r = 1}^t
            }
        }.
    \end{align*}
    Since \(\sigma\) only reorders the multiset of bucket projections, the median is unchanged:
    \begin{equation*}
        \med \lprp{
            \lbrb{
                \inp{P_{\sigma,\phi}^\top v}{\wt{Z}_{\sigma(r)}}
            }_{r = 1}^t
        }
        =
        \med \lprp{
            \lbrb{
                \inp{P_{\sigma,\phi}^\top v}{\wt{Z}_{r}}
            }_{r = 1}^t
        }.
    \end{equation*}
    Since \(P_{\sigma,\phi}\) is a permutation matrix \(\norm{P_{\sigma,\phi}^\top v}_q = \norm{v}_q\). Hence, after the change of variables \(u = P_{\sigma,\phi}^\top v\),
    \begin{equation*}
        \wt{F}(P_{\sigma,\phi} x)
        =
        \max_{u: \norm{u}_q = 1}
        \abs*{
            \inp{u}{x}
            -
            \med \lprp{
                \lbrb{
                    \inp{u}{\wt{Z}_{r}}
                }_{r = 1}^t
            }
        }
        =
        \wt{F}(x).
    \end{equation*}
\end{proof}

\begin{corollary}
    \label{cor:all_active_equal}
    There exists a minimizer \(x^\star\) of \(\wt{F}\) and a scalar \(a \in \mathbb{R}\) such that
    \begin{equation*}
        \forall j \in \bigcup_{r = 1}^t C_r: x^\star_j = a.
    \end{equation*}
\end{corollary}
\begin{proof}
    Start from the minimizer given by \cref{cor:bucketwise_constant}, so that
    \begin{equation*}
        \forall r \in [t], \; \forall j \in C_r: x_j = a_r
    \end{equation*}
    for some scalars \(a_1, \dots, a_t\).

    Let \(\mathcal{G}_{\mathrm{buck}}\) denote the finite set of all permutations \(P_{\sigma,\phi}\) from \cref{lem:bucket_swap_symmetry}. Define
    \begin{equation*}
        \bar{x} \coloneqq \frac{1}{|\mathcal{G}_{\mathrm{buck}}|} \sum_{P \in \mathcal{G}_{\mathrm{buck}}} P x.
    \end{equation*}
    By convexity of \(\wt{F}\) from \cref{obs:f_conv,lem:bucket_swap_symmetry},
    \begin{equation*}
        \wt{F}(\bar{x})
        \leq
        \frac{1}{|\mathcal{G}_{\mathrm{buck}}|} \sum_{P \in \mathcal{G}_{\mathrm{buck}}} \wt{F}(Px)
        =
        \wt{F}(x).
    \end{equation*}
    Hence \(\bar{x}\) is also a minimizer.

    The group \(\mathcal{G}_{\mathrm{buck}}\) acts transitively on the active coordinates \(\cup_{r = 1}^t C_r\): any active coordinate can be mapped to any other active coordinate by first permuting buckets and then using a bijection inside the target bucket. Therefore, \(\bar{x}\) must take the same value on all active coordinates, i.e.,
    \begin{equation*}
        \exists a \in \mathbb{R} \; \text{such that} \; \forall j \in \bigcup_{r = 1}^t C_r: \bar{x}_j = a.
    \end{equation*}
\end{proof}

\begin{lemma}
    \label{lem:inactive_zero}
    There exists a minimizer \(x^\star\) of \(\wt{F}\) such that
    \begin{equation*}
        \forall j \notin \bigcup_{r = 1}^t C_r: x^\star_j = 0.
    \end{equation*}
    Consequently, there exists a minimizer of the form
    \begin{equation*}
        x^\star = a \sum_{r = 1}^t \sum_{j \in C_r} e_j
    \end{equation*}
    for some scalar \(a \in \mathbb{R}\).
\end{lemma}
\begin{proof}
    Fix \(j \notin \cup_{r = 1}^t C_r\), and let \(D_j\) be the diagonal matrix defined by
    \begin{equation*}
        (D_j)_{k,k}
        =
        \begin{cases}
            -1 & \text{if } k = j, \\
            1 & \text{otherwise}.
        \end{cases}
    \end{equation*}
    Since every \(\wt{Z}_r\) is supported on \(C_r\), we have
    \begin{equation*}
        \forall r \in [t]: D_j \wt{Z}_r = \wt{Z}_r.
    \end{equation*}
    Hence, arguing exactly as in \cref{lem:within_bucket_symmetry},
    \begin{equation*}
        \forall x \in \mathbb{R}^d: \wt{F}(D_j x) = \wt{F}(x).
    \end{equation*}

    Let \(x\) be any minimizer of \(\wt{F}\), and define
    \begin{equation*}
        \bar{x} \coloneqq \frac{x + D_j x}{2}.
    \end{equation*}
    By convexity of \(\wt{F}\) from \cref{obs:f_conv},
    \begin{equation*}
        \wt{F}(\bar{x})
        \leq
        \frac{\wt{F}(x) + \wt{F}(D_j x)}{2}
        =
        \wt{F}(x).
    \end{equation*}
    Thus \(\bar{x}\) is also a minimizer. Moreover, \(\bar{x}\) agrees with \(x\) on every coordinate except \(j\), and
    \begin{equation*}
        \bar{x}_j = 0.
    \end{equation*}
    Repeating this argument for every \(j \notin \cup_{r = 1}^t C_r\) yields a minimizer supported on \(\cup_{r = 1}^t C_r\).

    Combining this with \cref{cor:all_active_equal} proves that there exists a minimizer of the form
    \begin{equation*}
        x^\star = a \sum_{r = 1}^t \sum_{j \in C_r} e_j = a \sum_{i = 1}^n e_{c(i)}.
    \end{equation*}
\end{proof}
At this point, by \cref{lem:inactive_zero}, it remains only to optimize over the scalar multiplying the active coordinates. Define
\begin{equation*}
    \rho \coloneqq \left\lfloor \frac t2\right\rfloor + 1,
    \qquad
    \eta \coloneqq t - \rho.
\end{equation*}
Let \(y_{(1)} \geq \cdots \geq y_{(t)}\) denote the non-increasing rearrangement of a multiset of scalars \(\bm{Y} = \{y_i\}_{i = 1}^t\). We use the convention
\begin{equation*}
    \med(\bm{Y}) =
    \begin{cases}
        y_{(k)} & \text{if } t = 2k - 1, \\
        \dfrac{y_{(k)} + y_{(k + 1)}}{2} & \text{if } t = 2k.
    \end{cases}
\end{equation*}
For \(a \in \mathbb{R}\), define
\begin{equation*}
    x^{(a)} \coloneqq a \sum_{i = 1}^n e_{c(i)} = a \sum_{r = 1}^t \sum_{j \in C_r} e_j.
\end{equation*}

\begin{lemma}
    \label{lem:bucket_sum_reduction}
    For every \(a \in \mathbb{R}\),
    \begin{equation*}
        \wt{F}(x^{(a)})
        =
        m^{1/p}
        \max_{u \in \mathbb{R}^t: \norm{u}_q \leq 1}
        \abs*{
            a \sum_{r = 1}^t u_r
            -
            \frac{\sqrt{d}}{m} \med \lprp{\lbrb{u_r}_{r = 1}^t}
        }.
    \end{equation*}
\end{lemma}
\begin{proof}
    For \(v \in \mathbb{R}^d\), define
    \begin{equation*}
        u_r \coloneqq m^{-1/p} \sum_{j \in C_r} v_j,
        \qquad r \in [t].
    \end{equation*}
    By Hölder,
    \begin{equation*}
        \abs*{u_r}
        =
        m^{-1/p} \abs*{\sum_{j \in C_r} v_j}
        \leq
        \norm{v_{C_r}}_q.
    \end{equation*}
    Hence,
    \begin{equation*}
        \sum_{r = 1}^t \abs*{u_r}^q
        \leq
        \sum_{r = 1}^t \norm{v_{C_r}}_q^q
        \leq
        \norm{v}_q^q.
    \end{equation*}
    In particular, every \(v\) with \(\norm{v}_q \leq 1\) produces a vector \(u\) with \(\norm{u}_q \leq 1\). Moreover,
    \begin{equation*}
        \sum_{j \in C_r} v_j = m^{1/p} u_r
        \qquad \text{for all } r \in [t].
    \end{equation*}
    Therefore,
    \begin{align*}
        \wt{F}(x^{(a)})
        &=
        \max_{v: \norm{v}_q = 1}
        \abs*{
            a \sum_{r = 1}^t \sum_{j \in C_r} v_j
            -
            \frac{\sqrt{d}}{m}
            \med \lprp{
                \lbrb{
                    \sum_{j \in C_r} v_j
                }_{r = 1}^t
            }
        } \\
        &=
        m^{1/p}
        \max_{v: \norm{v}_q = 1}
        \abs*{
            a \sum_{r = 1}^t u_r
            -
            \frac{\sqrt{d}}{m}
            \med \lprp{\lbrb{u_r}_{r = 1}^t}
        } \\
        &\leq
        m^{1/p}
        \max_{u: \norm{u}_q \leq 1}
        \abs*{
            a \sum_{r = 1}^t u_r
            -
            \frac{\sqrt{d}}{m}
            \med \lprp{\lbrb{u_r}_{r = 1}^t}
        }.
    \end{align*}
    For the reverse inequality, fix any \(u \in \mathbb{R}^t\) with \(\norm{u}_q \leq 1\), and define \(v \in \mathbb{R}^d\) by setting
    \begin{equation*}
        \forall r \in [t], \; \forall j \in C_r: v_j \coloneqq m^{-1 + 1/p} u_r,
    \end{equation*}
    and \(v_j \coloneqq 0\) outside \(\cup_{r = 1}^t C_r\). Then
    \begin{equation*}
        \sum_{j \in C_r} v_j = m^{1/p} u_r
        \qquad \text{for all } r \in [t],
    \end{equation*}
    and
    \begin{equation*}
        \norm{v}_q^q
        =
        \sum_{r = 1}^t \sum_{j \in C_r} m^{-q + q/p} \abs*{u_r}^q
        =
        \sum_{r = 1}^t m^{(1 - q) + q/p} \abs*{u_r}^q
        =
        \sum_{r = 1}^t \abs*{u_r}^q
        =
        \norm{u}_q^q
        \leq 1,
    \end{equation*}
    since \((1 - q) = -q / p\) as $1/p + 1/q = 1$. Plugging $v$ into the definition of $\wt{F}$ now yields the reverse inequality, concluding the proof of the lemma.
\end{proof}

\begin{lemma}
    \label{lem:exact_scalar_objective}
    For every \(a \in \mathbb{R}\),
    \begin{equation*}
        \wt F(x^{(a)})
        =
        m^{1/p}
        \left(
            \rho
            \abs*{\frac{\sqrt d}{\rho m}-a}^p
            +
            \eta |a|^p
        \right)^{1/p}.
    \end{equation*}
\end{lemma}

\begin{proof}
    Let
    \[
        \lambda \coloneqq \frac{\sqrt d}{m}.
    \]
    By \cref{lem:bucket_sum_reduction}, it suffices to compute
    \[
        \max_{\norm{u}_q\leq 1}
        \abs*{
            a\sum_{r=1}^t u_r
            -
            \lambda \med\lprp{\lbrb{u_r}_{r=1}^t}
        } .
    \]
    Since \(\med(-u)=-\med(u)\), this maximum
    equals
    \[
        \max_{\norm{u}_q\leq 1}
        \left\{
            \lambda \med\lprp{\lbrb{u_r}_{r=1}^t}
            -
            a\sum_{r=1}^t u_r
        \right\}.
    \]
    Let
    \[
        u_{(1)}\geq u_{(2)}\geq\cdots\geq u_{(t)}
    \]
    be the non-increasing rearrangement of \(u\). Now, observe that:
    \[
        \med(u)
        \leq
        \frac1\rho \sum_{r=1}^{\rho} u_{(r)} .
    \]
    Therefore,
    \begin{equation*}
        \lambda \med(u)-a\sum_{r=1}^t u_r
        \leq
        \frac{\lambda}{\rho}\sum_{r=1}^{\rho}u_{(r)}
        -
        a\sum_{r=1}^t u_{(r)} =
        \left(\frac{\lambda}{\rho}-a\right)
        \sum_{r=1}^{\rho}u_{(r)}
        -
        a\sum_{r=\rho+1}^t u_{(r)} .
    \end{equation*}
    By Hölder,
    \[
        \lambda \med(u)-a\sum_{r=1}^t u_r
        \leq
        \left(
            \rho
            \abs*{\frac{\lambda}{\rho}-a}^p
            +
            \eta |a|^p
        \right)^{1/p}
        \norm{u}_q .
    \]
    Since \(\norm{u}_q\leq 1\), this establishes an upper bound.

    For the lower bound, set
    \[
        A \coloneqq \frac{\lambda}{\rho}-a,
        \qquad
        B \coloneqq -a, \qquad
        \alpha
        \coloneqq
        \left(
            \rho |A|^p+\eta |B|^p
        \right)^{-1/q} \text{, and }
        J(z)\coloneqq \operatorname{sgn}(z)|z|^{p-1}.
    \]
    Now, define \(u\in\mathbb R^t\) by
    \[
        u_1=\cdots=u_\rho=\alpha J(A),
        \qquad
        u_{\rho+1}=\cdots=u_t=\alpha J(B).
    \]
    Since $A - B > 0$, and \(J\) is increasing, this vector is non-increasing. Also, observe
    \[
        \norm{u}_q^q
        =
        \alpha^q
        \left(
            \rho |A|^p+\eta |B|^p
        \right)
        =
        1 \text{ and } \med(u)=\alpha J(A).
    \]
    Hence
    \begin{align*}
        \lambda \med(u)-a\sum_{r=1}^t u_r
        &=
        \lambda \alpha J(A)
        -
        a\left(
            \rho \alpha J(A)+\eta \alpha J(B)
        \right) =
        \alpha
        \left(
            \rho A J(A)+\eta B J(B)
        \right) \\
        &=
        \alpha
        \left(
            \rho |A|^p+\eta |B|^p
        \right) =
        \left(
            \rho |A|^p+\eta |B|^p
        \right)^{1/p}.
    \end{align*}
    matching the upper bound and proving the lemma.
\end{proof}

\begin{corollary}
    \label{cor:optimal_scalar}
    The optimal scalar is
    \[
        a^\star
        =
        \frac{\sqrt d}{\rho m}
        \cdot
        \frac{1}{
            1+\left(\eta/\rho\right)^{1/(p-1)}
        }.
    \]
    Consequently, a minimizer of \(\wt F\) is
    \[
        x^\star
        =
        a^\star \sum_{i=1}^n e_{c(i)},
    \]
    and hence, by \cref{obs:sym_f}, a minimizer of \(F\) is
    \[
        \wh{\mu}_{\mrm{proj}} (\bm{X})
        =
        S x^\star
        =
        a^\star \sum_{i=1}^n s_i e_{c(i)} .
    \]
\end{corollary}
\begin{proof}
    Observing that $\wt{F}(x^{(a)})$ is convex in $a$ and computing its derivative, the corollary follows.
\end{proof}
\begin{corollary}
    \label{cor:norm_of_minimizer}
    The resulting estimator has
    \begin{equation*}
        \norm{\wh{\mu}_{\mrm{proj}} (\bm{X})}_p = \norm{x^\star}_p
        =
        a^\star n^{1/p}
        \geq
        \frac{\sqrt{d}}{2 n^{1 - 1/p}}.
    \end{equation*}
\end{corollary}
\begin{proof}
    Since \(x^\star\) has exactly \(n\) nonzero coordinates, each of magnitude \(a^\star\),
    \begin{equation*}
        \norm{x^\star}_p = a^\star n^{1/p}.
    \end{equation*}
    The corollary follows from \cref{cor:optimal_scalar}.
\end{proof}
\begin{proof}[Proof of \cref{thm:proj_sub_optimal}]
    The theorem follows from \cref{cor:norm_of_minimizer,lem:unique_coords}.
\end{proof}

\section{Proof of \cref{thm:mom_l2_lb_main}}
\label{sec:suboptimalcoordwise}
In this section, we show that \cref{alg:MOM} cannot be optimal for the $\ell_2$ norm.
The optimal rate for heavy-tailed mean estimation for isotropic data in this setting is known to be of order
\begin{equation}\label{eq:optimal_l2}
\sqrt{\frac{d}{n}} + \sqrt{\frac{1}{n} \log\left(\frac{1}{\delta}\right)}.
\end{equation}

We prove the following lower bound which demonstrates that the coordinate-wise median-of-means strategy is sub-optimal when the failure probability is small.
\MoMSuboptimal*
This result follows from combining two different lower bounds: the first demonstrates a large error when \cref{alg:MOM} is run with many buckets, while the second shows a large error when it is run with few buckets.

\begin{lemma}\label{lem:lowerbound_constant_buckets}
There is a universal constants $C > 0$ such that the following holds.
For any $\delta \in (0, 1/3)$ and $ n, d, t \in\mathbb{N}$ with $ t\leq n $, there exists a distribution $ P_1 $ over $\mathbb{R}^{d} $ with mean $0$ and covariance $I_d$ such that if \cref{alg:MOM} is given the number of buckets $t$ and the sample $X_1, \dots, X_n \sim P_1$ as input, it returns an estimate $\hat \mu$ satisfying
\begin{align*}
    \Pr_{X_1, \dots, X_n \sim P_1} \left[
    \left\|\hat\mu\right\|_{2}
    \geq
    C\sqrt{\frac{dt}{n}}
    \right]
    > \delta .
\end{align*}
\end{lemma}

\begin{proof}
Set $p = 1-\frac{t}{8n}$ and $\ell=\frac{1}{\sqrt{p(1-p)}}$ and define  random variable $Z \sim Q$ by
\[
    \Pr_{Z \sim Q}[Z=(1-p)\ell]=p,
    \qquad
    \Pr_{Z \sim Q}[Z=-p\ell]=1-p.
\]
Notice that $\E_{Z \sim Q}[Z] = 0$ and $Z$ has variance $1$.
Let $P_1$ be the product distribution on $\mathbb R^d$ whose coordinates are
independent copies of $Z$. Clearly $P_1$ has mean $0$ and covariance $I_d$.

Let $X_1,\ldots,X_n$ be an i.i.d.\ sample from $P_1$ and run \cref{alg:MOM} using $t$
buckets. Denote the bucket sizes by $m_1,\ldots,m_t$, so that
\[
    m_j\leq \left\lceil \frac nt\right\rceil \leq \frac nt+1 .
\]
Fix a coordinate $k\in[d]$ and a bucket $j\in[t]$. The empirical mean of the
$k$th coordinate in bucket $j$ equals $(1-p)\ell$ whenever all samples in that
bucket have $k$th coordinate equal to $(1-p)\ell$. This event has probability
\[
    p^{m_j}
    =
    \left(1-\frac{t}{8n}\right)^{m_j}
    \geq
    1-\frac{m_jt}{8n}
    \geq
    1-\left(\frac nt+1\right)\frac{t}{8n}
    \geq
    \frac34,
\]
where we used the inequality $(1+x)^y \ge 1+yx$ which holds for $x \ge -1$ and $ y \ge 1$, and the fact that $t \le n$.
For the fixed coordinate $k \in [d]$, let $Y_k$ be the number of buckets whose
coordinate-$k$ empirical mean is equal to $(1-p)\ell$. Then $Y_k$ is a sum of
independent Bernoulli random variables with
\[
    \E Y_k\geq \frac{3t}{4}.
\]
We can apply Markov's inequality to the random variable $t - Y_k$ to get
\[
   \Pr_{X_1, \dots, X_n \sim P_1}\left[Y_k>\frac t2\right] = 1- \Pr_{X_1, \dots, X_n \sim P_1}\left[t-Y_k\geq \frac t2\right] \geq 1 - \frac{\frac{t}{4}}{\frac{t}{2}}= \frac{1}{2}.
\]

On the event $\{Y_k>t/2\}$, strictly more than half of the bucket means in
coordinate $k$ are equal to $(1-p)\ell$, and therefore the coordinate-wise
median $M_k$ is $(1-p)\ell$. Now define
\[
    W=\sum_{k=1}^d \mathbf{1}\{M_k=(1-p)\ell\}.
\]
Since the coordinates of $P_1$ are independent, the events
$\{M_k=(1-p)\ell\}$ are independent over $k\in[d]$. Moreover since,
\[
    \E W\geq \frac{d}{2},
\]

we can invoke Markov's inequality on the random variable $d-W$ to conclude that
\[
    \Pr_{X_1, \dots, X_n \sim P_1}\left[W\geq \frac{d}{4}\right] = 1-\Pr_{X_1, \dots, X_n \sim P_1}\left[d-W\geq \frac{3d}{4}\right] \ge 1- \frac{\frac{d}{2}}{\frac{3d}{4}} = \frac{1}{3}.
\]
On this event, we have that
\[
    \|\hat\mu\|_2
    \geq
    \sqrt{\frac{d}{4}}(1-p)\ell = \sqrt{\frac{dt}{4(8n-t)}} \ge C\sqrt{\frac{dt}{n}},
\]
for some universal constant $C > 0$, where we have used the assumption $ t\leq n $.
By setting $c = \frac{1}{3}$ it holds for every $\delta\in(0,c)$ that
\[
    \Pr_{X_1,\ldots,X_n\sim P_1}
    \left[
        \|\hat\mu\|_2
        \geq
        C\sqrt{\frac{dt}{n}}
    \right]
    >
    \delta,
\]
as claimed. This completes the proof.
\end{proof}

\begin{lemma}
\label{lem:radial_rademacher_mom_small_delta}
There exist universal constants $c_1, c_2$ such that the following holds.
Let $n,d,t\in\mathbb N$ with $t\le n$, and
$\delta\in(0, c_1)$ satisfying
\[
    t\le \frac{\log(1/\delta)}{\log\log(1/\delta)}.
\]
There exists a distribution $P$ over $\mb R^d$ satisfying $\E_P X=0$ and $\Cov_P(X)=I_d$, such that the following holds.
If \cref{alg:MOM} is run with $t$ buckets on
$X_1,\ldots,X_n\sim P$, its output $\hat\mu$ satisfies
\[
    \Pr_{X_1,\ldots,X_n\sim P}
    \left[
        \norm{\hat\mu}_2
        \ge
        c_2\sqrt{\frac{td\log(1/\delta)}{n}}
    \right]
    >
    \delta .
\]
\end{lemma}

\begin{proof}
Our hard distribution will be of the form
\[
    X=R\varepsilon,
    \qquad
    \varepsilon=(\varepsilon_1,\ldots,\varepsilon_d),
\]
where $\varepsilon_1,\ldots,\varepsilon_d$ are independent Rademacher random variables and $R$ is an independent scalar random variable to be specified subsequently. Now, let
\[
    L=\log(1/\delta),
    \qquad
    m=\left\lceil \frac nt\right\rceil,
    \qquad
    a=\frac{L}{8e^2}.
\]
Note that $a \ge 2$ for sufficiently small $c_1$. Now, define $R$ as follows
\[
    \Pr\left[R=\sqrt{ma}\right]=\frac{1}{ma},
    \qquad
    \Pr[R=0]=1-\frac{1}{ma},
\]
and let $X=R\varepsilon$. Then, observe that
\[
    \E X=0,
    \qquad
    \Cov(X)=\E[R^2]\E[\varepsilon\varepsilon^\top]
    =
    ma\cdot \frac{1}{ma} I_d
    =
    I_d.
\]

Let $m_1,\ldots,m_t$ be the bucket sizes produced by \cref{alg:MOM}. Since
$t\le n$,
\[
    m_i\in
    \left\{
        \left\lfloor\frac nt\right\rfloor,
        \left\lceil\frac nt\right\rceil
    \right\},
    \qquad
    m_i\ge \frac{m}{2}.
\]
For bucket $i$, let $E_i$ be the event that exactly one sample in the bucket has nonzero multiplier $R$. Then
\[
    \Pr[E_i]
    =
    m_i \cdot \frac{1}{ma}
    \cdot \left(1-\frac{1}{ma}\right)^{m_i-1}
    \ge
    \frac{1}{8a}.
\]
Observing that $E_1,\ldots,E_t$ are independent, we get for
$E=\cap_{i=1}^t E_i$,
\begin{equation}
\label{eq:radial_rademacher_good_event}
    \Pr[E]\ge \left(\frac{1}{8a}\right)^t .
\end{equation}

We condition on $E$ for the remainder of the proof. Note that each bucket contains exactly one nonzero sample, and the corresponding Rademacher vectors remain independent. Hence the $i$th bucket mean is distributed as
\[
    \hat\mu^{(i)}
    =
    \frac{\sqrt{ma}}{m_i}\varepsilon^{(i)}
    =
    s_i\varepsilon^{(i)},
    \qquad
    i=1,\ldots,t,
\]
where $\varepsilon^{(1)},\ldots,\varepsilon^{(t)}$ are independent copies of $\varepsilon$ and
\[
    s_i\ge s_0\coloneqq \sqrt{\frac{a}{m}}.
\]
Thus, conditioned on $E$,
\[
    \hat\mu_j
    =
    \mathrm{Median}
    \left(
        s_1\varepsilon^{(1)}_j,\ldots,s_t\varepsilon^{(t)}_j
    \right),
    \qquad j\in[d].
\]

We use the following elementary fact.
\begin{claim}
    \label{clm:rademacher_median}
    For every $s_i\ge s_0$,
    \[
        \Pr\left[
            \left|
            \mathrm{Median}(s_1\varepsilon_1,\ldots,s_t\varepsilon_t)
            \right|
            \ge s_0
        \right]
        \ge 1/2,
    \]
    where $\varepsilon_1,\ldots,\varepsilon_t$ are independent Rademacher random variables.
\end{claim}
\begin{proof}
    If $t$ is odd, one of the two signs appears strictly more than $t/2$ times. The median then has that sign and absolute value at least $s_0$. If $t$ is even, the same conclusion holds unless exactly half of the signs are positive and half are negative. Therefore
    if $B\sim \mrm{Binom}(t,1/2)$, then
    \[
        \Pr\left[
            \left|
            \mathrm{Median}(s_1\varepsilon_1,\ldots,s_t\varepsilon_t)
            \right|
            \ge s_0
        \right]
        \ge
        1-\ind\{t \text{ even}\}
        \Pr\left[
            B=\frac t2
        \right]
        \ge \frac12 .
    \]
\end{proof}
Let $Z=\sum_{j=1}^d \ind\left\{|\hat\mu_j|\ge s_0\right\}$. Conditioned on $E$, the summands defining $Z$ are independent and have expectations at least $1/2$ by \cref{clm:rademacher_median}. Hence a Chernoff bound gives
\[
    \Pr\left[
        Z\ge \frac d4
        \,\middle|\,
        E
    \right]
    \ge c.
\]
On this event, $\norm{\hat\mu}_2\ge s_0\sqrt{d}/2$. Therefore
\[
    \Pr\left[
        \norm{\hat\mu}_2
        \ge
        c\,s_0\sqrt{d}
        \,\middle|\,
        E
    \right]
    \ge c.
\]
Since $m\le 2n/t$,
\[
    s_0\sqrt d
    =
    \sqrt{\frac{a}{m}}\sqrt d
    \ge
    c\sqrt{\frac{atd}{n}}.
\]
Therefore, using \cref{eq:radial_rademacher_good_event},
\[
    \Pr\left[
        \norm{\hat\mu}_2
        \ge
        c\sqrt{\frac{atd}{n}}
    \right]
    \ge
    c\left(\frac{1}{8a}\right)^t
    =
    c\left(\frac{e^2}{L}\right)^t .
\]
By taking $ c_{1} $ sufficently small such that $ \log L-2\geq0 $  and since $t\le L/\log L$,
\[
    \log\left[
        c\left(\frac{e^2}{L}\right)^t
    \right]
    \ge
    \log c
    -
    \frac{L}{\log L}(\log L-2)
    =
    -L+\frac{2L}{\log L}+\log c.
\]
Taking $c_1$ sufficiently small makes the last quantity larger than $-L$.
Hence
\[
    c\left(\frac{e^2}{L}\right)^t
    >
    e^{-L}
    =
    \delta.
\]
Finally, because $a=L/(8e^2)$,
\[
    c\sqrt{\frac{atd}{n}}
    =
    c'\sqrt{\frac{tdL}{n}}
    =
    c'\sqrt{\frac{td\log(1/\delta)}{n}},
\]
for another universal constant $c'>0$. Renaming constants establishes the lemma.
\end{proof}

\begin{proof}[Proof of \cref{thm:mom_l2_lb_main}]
    The result follows from \cref{lem:lowerbound_constant_buckets,lem:radial_rademacher_mom_small_delta} for the regimes:
    \begin{equation*}
        t \geq \frac{\log (1 / \delta)}{\log \log (1 / \delta)} \text{ and } t \leq \frac{\log (1 / \delta)}{\log \log (1 / \delta)}
    \end{equation*}
    respectively.
\end{proof}
\end{document}